\documentclass[12pt]{article}
\usepackage{amsfonts}
\usepackage{amssymb}
\usepackage{amscd, color}
\usepackage{hyperref}
\usepackage{amsmath}
\usepackage{mathabx}

\newcommand{\pref}[1]{\ref{#1} page \pageref{#1} }

\def\boxit#1#2{\setbox1=\hbox{\kern#1{#2}\kern#1}%
\dimen1=\ht1 \advance\dimen1 by #1 \dimen2=\dp1 \advance\dimen2 by #1
\setbox1=\hbox{\vrule height\dimen1 depth\dimen2\box1\vrule}%
\setbox1=\vbox{\hrule\box1\hrule}%
\advance\dimen1 by .4pt \ht1=\dimen1
\advance\dimen2 by .4pt \dp1=\dimen2 \box1\relax}

\newcommand{\BEQ}{\begin{equation}}
\newcommand{\EEQ}{\end{equation}}
\newtheorem{Th}{Theorem}

\newtheorem{Lem}{Lemma}
\newtheorem{Prop}{Proposition}
\newtheorem{Cor}{Corollary}

\newtheorem{Ex}{Example}

\newtheorem{Def}{Definition}
\newtheorem{Rem}{Remark}

\def\nn{\nonumber}
\def\one#1{#1^{\raise5pt\hbox{$\scriptstyle\!\!\!\!1$}}\,{}}
\def\two#1{#1^{\raise5pt\hbox{$\scriptstyle\!\!\!\!2$}}\,{}}
\def\oneM#1{#1^{\raise5pt\hbox{$\scriptstyle\!\!\!\!m$}}\,{}}
\def\oneK#1{#1^{\raise5pt\hbox{$\scriptstyle\!\!\!\!k$}}\,{}}

\def\bea{\begin{eqnarray}}
\def\eea{\end{eqnarray}}

\newcommand{\CC}{{{\mathbb{C}}}}
\newcommand{\kk}{k}
\newcommand{\kki}{\kk[i]}

\def\RR{{\mathbb{ R}}}

\def\mpsi{ \psi^M }

\def\MMs{Manin\ matrices}
\def\MM{Manin\ matrix}
\def\BX{$\Box$}
\def\PRF{ {\bf Proof.} }

\def\DD{D} %
\newtheorem{Notn}{Notation}

\begin{document}


\begin{titlepage}
\hfill ITEP-TH-47/11

\vskip 1cm

\centerline{\Large Decomplexification of the Capelli identities and holomorphic factorization}

\vskip 1cm
\centerline{Preliminary version. Comments are welcome at \href{http://alexchervov.wordpress.com/2012/03/08/}{http://alexchervov.wordpress.com/2012/03/08/} }

\vskip 1cm
 \centerline{A. Chervov \footnote{E-mail:
al.chervov@gmail.com}
}


\vskip 1.0cm
\centerline{\large \bf  Abstract}
\vskip 1cm

The Capelli identities claim $det(A)det(B) = det(AB+correction)$ for certain matrices with noncommutative entries. They have applications in representation theory and integrable systems. 
Several examples are known, but the  general theory is not yet found.  
We propose new examples of these identities,
constructed according to the following principle.  For several known identities for $n$ by $n$
matrices we construct new identity for $2n$ by $2n$ matrices where each element $z$ of the original matrix is substituted by 2x2 matrix of the form $[real(z) ~imag(z);~ -imag(z)~ real(z)]$,
i.e. we view the original identity as complex valued, while the new identity is its real form (decomplexification). 
It appears  that "decomplexification" affects non-trivially  the "correction term".  It becomes tridiagonal matrix, in contrast to the diagonal in the classical case.
The key result is an extension to the non-commutative setting of the fact that the determinant of the decomplexified matrix is equal to the square module of the determinant of the original matrix (in non-commutative setting the corrections are necessary). The decomplexified Capelli's identities are 
 corollaries of this fact and standard Capelli identities. 
We also discuss  analogues of the Cayley identity;  observe that the Capelli determinant
coincides with the Duflo image of the standard determinant; give  short proof of the Cayley identity
via  Harish-Chandra's radial part calculation.

The main motivation for us is a recent paper by An Huang  
(\href{http://arxiv.org/abs/1102.2657}{arXiv:1102.2657}). 
From our viewpoint his result is a "dequaternionification" 
of $1$ by $1$ Capelli identity. Apparently it
can be extended to $n$ by $n$ case, but our approach should be somehow modified for this. 

The 
paper aims to be accessible and interesting not only for experts. It gives brief review
of Capelli identities, applications, their relation with Wick quantization, Duflo map, some open issues, etc.

\vskip 1cm



{\bf Mathematics Subject Classification (MSC) codes: } 
15B99, 05A19, 15A15, 17B35, 20G05
15B99 
05A19 
15A15 
17B35  
20G05 


\end{titlepage}

\tableofcontents


\section{Introduction}
Introduction contains 4 subsections 1) short summary of main messages of the paper 2) description of new results 3) overview of  applications, related questions and surrounding results 4) organization of paper,
ideas of the proofs, informal comments.

\subsection{Short summary}

The Capelli identities claim $det(A)det(B) = det(AB+correction)$ for certain matrices with noncommutative entries.
Here we propose new results of that kind. 

The identities are explicit and somewhat nice formulas which any undergraduate can understand.
So we try to make paper understandable and interesting not only for experts. We include some background
material, discussion of the broader context, etc. 
Let us start with a short list of highlights  which might be of interest - they contain  overview of results, context, etc.

1) Main results. (See subsection \pref{ResultSubSect} for the list of results). \\
Given a Capelli identity $det(A)det(B) = det(AB+correction)$  for $n$ by $n$ 
matrices, we construct new identity for $2n$ by $2n$ matrices, where each element $z$ of the original matrices is substituted by the 2x2 matrix of the form $[real(z) ~imag(z);~ -imag(z)~ real(z)]$.
I.e. we view the original identity as complex valued, and the new identity as is its real form (decomplexification). 
We generalize the standard Capelli identity, its symmetric and anti-symmetric analogues, and recently found quite general
identities from \cite{CSS08}.
We also propose an analogue for the factorization theorem: $det(M^\RR+cor_R)= det(M+cor)det(\bar M+\widebar{cor})$. 
This result is our key tool to approach Capelli identities.
We also discuss  analogues (decomplexifications, dequaternionification) of the Cayley identity
($det(d/dx_{ij} )det(x_{ij} )^s = s(s+1)...(s+n-1) det(x_{ij} )^{s-1}$), 
although they appear to be rather direct.  Some open problems are discussed.

2) General context. Determinants and linear algebra for {\bf good} matrices with non-commuting elements.
(See subsection \pref{StoryAboutSubSect} for further comments). \\ 
It is rather clear that there is no reasonable notion of the determinant with values in the ground ring for  matrices with {\bf generic} non-commutative entries and there is no hope to get  direct analogues of linear algebra theorems. 
 However the situation is different if we consider only certain specific classes of matrices
with non-commutative entries. 
The best known 
example is quantum groups and $q$-determinant.
Capelli matrix and Capelli determinant; super-matrices and Berezinian and also  Manin matrices \cite{CFR09}   are other examples. 
So we can pose a question: {\bf what matrices with non-commuting elements are  {
"good" }  } ? 
(i.e. "good" - means:  they have the determinant and linear algebra theorems (e.g. Cayley-Hamilton, Cramer, Cauchy-Binet, etc.)  holds true for them (this has been done for Manin matrices in 
\cite{CFR09})).

3) Why matrices with non-commutative elements are interesting ? (See subsection \pref{sectWhy}).\\
A) Lax matrices of quantum integrable systems are like this, and we want to construct
quantum commuting integrals of motion by analogues of classical formulas $det(L(z))$, $Tr L^k(z)$.
B) In representation theory interesting elements in universal enveloping algebras can be
constructed by similar formulas (see section \ref{sectWhy}, page \pageref{sectWhy}).
C) We may dream about real world applications of such matrices, since non-commuting
elements can be matrices themselves - so we can dream to find some real world application with
ordinary matrices
which after splitting into blocks become examples of ``good'' matrices over the matrix ring.
Numerical algorithms for such matrices can be modified to work faster - this would be of value
for real world applications (if such exist).

4) Capelli identities mean - Wick ordering is great ! (See subsection \pref{WickIntrSect}). \\
Wick ordering in quantization theory is a linear map which sends a commutative monomial $z^ap^b$ to a
 non-commutative monomial $z^a \partial_z^b$ i.e. derivatives stand on the right. 
Capelli identities can be reformulated as: 
\bea
F_{quantum} (Wick(ZP^t)) = Wick ( F_{classical} (ZP^t)),
\eea
where $F$ can be determinant, permanent, etc., $F_{quantum}$ is certain quantum analogue of it 
(e.g. Capelli determinant), $Z$ is a matrix with elements $z_{ij}$ and $P$ with $p_{ij}$.
In general guessing ``correct quantum analogues'' is a non-trivial art. 
So identities mean that  in certain situations we can get the correct quantum analogues
from the classical ones just taking Wick ordering of the corresponding classical expressions.
Which is very easy. But the proof that such recipe works 
cannot be achieved easily  - only via the Capelli identities.

5) Capelli determinant = Duflo (determinant). (See subsection \pref{CapDufloIntrSect}). \\
Duflo map is an amazing and powerful result working for any Lie algebra. Moreover M. Kontsevich \cite{Kon97} generalized it to quantization of any Poisson manifold. It allows to construct the central elements of non-commutative quantum algebras from the 
Poisson centers of corresponding commutative algebras. We observe (and it seems that this observation is new)
that the Capelli determinant is the Duflo image for $gl_n$ of the simple determinant of matrix with commutative elements.

6) Cayley identity via  Harish-Chandra's radial part 
calculation. (See appendix \pref{AppendSect})\\
In appendix we give  yet another very short proof of the Cayley identity via Harish-Chandra's radial part 
calculation of invariant differential operators. 
We put it here mainly for pedagogical reasons - it seems to us a nice exercise:
a short proof of a concrete result by means of general technique which recent days become
quite famous due to P. Etingof and V. Ginzburg \cite{EG00}, who
put it into a more general context and related to various topics
such as Calogero-Moser integrable hamiltonian system, Cherednik algebras, hamiltonian reduction, n!-conjecture,  etc.

\vskip .5cm

We try to make paper "modular" and "iterative". Modular means that "modules" are somewhat independent (i.e. paper should be readable from a random place, so one can start looking  directly at the place of interest, therefore each section starts with a brief overview; 
each subsection starts with brief formulation what going to be done. 
The presentation goes in 3 iterations: 
the short summary above is the first one,
Introduction being an extended summary and  is the second iteration, 
and the main text is the third one.
Any feedback on the paper is welcome by e-mail or   \href{http://alexchervov.wordpress.com/2012/03/08/}{http://alexchervov.wordpress.com/2012/03/08/}.


\subsection{Main results}
Here we describe main results - we have 3 subsections  with 1) general overview 
2) results devoted to Capelli identities (known results are recalled) 
3) results about factorization of the determinant $det(M^\RR+cor_R)= detr(M+cor)det(\bar M+\widebar{cor})$.

\subsubsection{What the story is about ? \label{StoryAboutSubSect}}
Let us briefly repeat what are the results about and what is the general context they fit (open problem
of understanding what are "good matrices" with non-commutative elements).

The Capelli identities state that $det(A)det(B) = det(AB+correction)$
\footnote{the original works were restricted to determinants, but later generalizations
to permanetns, immanants, pfaffians and traces were proposed, but we were not able to extend our results to cover these cases, so we restrict our discussion to determinants}  
for certain matrices with noncommutative entries. Such identities have applications in
representation theory and integrable systems (see subsections \ref{WickIntrSect}, \ref{sectWhy} for a brief comments on applications and \ref{OverCapSect} for brief overview of related works). 
In this paper we propose some new examples of these identities,
constructed by the following principle.  For several known identities for $n$ by $n$
matrices we construct new identity for $2n$ by $2n$ matrices, where each element $z$ of the original matrix is substituted by 2x2 matrix of the form $
\left( \begin{array}{ccc} real(z) & imag(z) \\ -imag(z)  & real(z) \end{array} \right).
$
I.e. we consider the original identity as complex valued, while the new identity is its real form ("decomplexification"). In more invariant terms, we consider a linear operator $M: \CC^n \to \CC^n$
and "decomplexified" matrix  $M^\RR$ 
 is the matrix of same operator considered as
acting: $\CC^n=\RR^{2n} \to \CC^n=\RR^{2n} $.
Somewhat non-obvious part of the work is to understand how the "correction"
term changes under "decomplexification" and proving new identities
(previous approaches do not seem to work).
The main motivation for us is a recent paper by   
\href{http://arxiv.org/abs/1102.2657}{An Huang} \cite{AH11}. 
From our viewpoint his result is a "dequaternionification" 
of $1$ by $1$ Capelli identity. Apparently it
can be extended to $n$ by $n$ case, but this is not achieved yet.

{\bf Broader conte$\frac{x}{s}$t.} Before going into details of this paper let us mention the following.
It is rather clear that there is no reasonable notion of the determinant with values in the ground ring for  matrices with
{\bf generic} non-commutative entries and there is no hope to get  direct analogues of linear algebra theorems.
\footnote{let us note that the theory of quasi-determinants  \cite{GGRW02}
deals with generic matrices, however quasi-determinants are not direct analogues
of the determinant, there is {\bf no unique} quasi-determinant, but $n^2$ of them
which are inverses to the elements of the inverse matrix. 
Quasi-determinantal linear algebra theorems sometimes quite far from the usual linear algebra theorems - but probably this is the best one can get in the framework
of generic matrices.}
 However the situation changes if we consider certain {\bf specific} matrices
with non-commutative entries, matrices whose elements satisfy certain
commutation relations. Then in certain examples 
it is possible to define "the determinant" and prove many linear algebra
theorems. The most famous example are quantum groups and $q$-determinant.
Capelli matrix and determinant and also  Manin matrices \cite{CFR09}   are other examples of this phenomena.
There are  other examples related to "Lax matrices" of quantum integrable
systems where similar constructions are known or expected.
However quite naive questions: {\bf for what matrices with non-commutative
entries there exists "the determinant"? if it exists, how should it be defined? } remain unanswered. 
Let us stress that the explicit formulas for the determinants depend on the commutation
relations of  matrix elements: for quantum groups - q-determinant, for super-matrices - Berezinian, in Capelli case - certain correction is required: det(E+correction).
The present paper 
provides new examples of Capelli type matrices and identities.
 It might be helpful for  building  the general theory.





\subsubsection{Decomplexified Capelli identities \label{ResultSubSect} }

Let us first recall known results, later formulate new ones.
Everywhere we denote by $D^t$ transpose to matrix $D$;  as ${\det}^{\mathrm{col}}$ - column-determinant of a matrix with not necessarily commuting elements. It is defined by the same rule as the standard one, with the convention that
elements from first column are at the first position in the products, from second column are at the second position, etc. See definition \pref{coldetdef} for more details.

Consider the polynomial algebra $\CC[z_{ij}]$, and the  $n$ by $n$   matrices:
\bea \label{ZFnot1}
Z = \left( \begin{array}{ccc}
z_{11} & ... & z_{1n} \\
... & ... & ... \\
z_{n1} & ... & z_{nn} \\
\end{array} \right), ~~~
\DD = \left( \begin{array}{ccc}
\frac{\partial}{\partial z_{11}}  & ... & \frac{\partial}{\partial z_{1n}} \\
... & ... & ... \\
\frac{\partial}{\partial z_{n1}}  & ... & \frac{\partial}{\partial z_{nn}} \\
\end{array} \right). ~~~
\eea

{\bf Known theorems} {
{
\it 

(\href{http://gdz.sub.uni-goettingen.de/index.php?id=11&no_cache=1&IDDOC=26783&IDDOC=26783&branch=&L=1}
{A. Capelli}
  \cite{Ca87})
\bea
{\det}^{\mathrm{col}} (Z \DD^t + \mathrm{diag}(n-1,n-2,...,1,0))= 
{\det} ( Z) {\det} (\DD^t).
\eea

 H.W. Turnbull \cite{Tu48} proved that the same formula holds true
 if $Z$ is the symmetric  matrix with $z_{ij}$ $i\le j$ on both positions $(ij)$ and $(ji)$ and respectively $\DD$ is the symmetric matrix with  $\partial_{z_{ij}}$ out of the diagonal and $2 \partial_{z_{ii}}$ on the diagonal.

R. Howe, T. Umeda \cite{HU91}  and  B. Kostant, S. Sahi \cite{KS91} \footnote{the explicit formulation which we use can be found in \cite{FZ93} }   covered the antisymmetric case,
i.e. if $n$ is even,  $Z$ is the anti-symmetric  matrix with $z_{ij}$ $i< j$ at $(ij)$ and $-z_{ij}$ at $(ji)$, and, respectively, $\DD$ is  corresponding anti-symmetric matrix,   then the same  identity holds.
}

Consider the polynomial algebra $\CC[x_{ij},y_{ij}]$, and the $2n$ by $2n$  matrices:
\bea
Z^{\RR} = \left( \begin{array}{ccccc}
x_{11} & y_{11} & ... & x_{1n} & y_{1n} \\
-y_{11} & x_{11} & ... & -y_{1n} & x_{1n} \\
... & ... & ... \\
x_{n1} & y_{n1} & ... & x_{nn} & y_{nn} \\
-y_{n1} & x_{n1} & ... & -y_{nn} & x_{nn} \\
\end{array} \right), ~~~
\DD^{\RR} = (\frac{1}{2}) \left( \begin{array}{ccccc}
\frac{\partial}{\partial x_{11}}  & -\frac{\partial}{\partial y_{11}}  & ...  & \frac{\partial}{\partial x_{n1}} & -\frac{\partial}{\partial y_{n1}} \\
\frac{\partial}{\partial y_{11}}  & \frac{\partial}{\partial x_{11}}  & ...  & \frac{\partial}{\partial y_{n1}} & \frac{\partial}{\partial x_{n1}} \\
... & ... & ... \\
\frac{\partial}{\partial x_{n1}}  & -\frac{\partial}{\partial y_{n1}}  & ...  & \frac{\partial}{\partial x_{nn}} & - \frac{\partial}{\partial y_{nn}} \\
\frac{\partial}{\partial y_{n1}}  & \frac{\partial}{\partial x_{n1}}  & ...  & \frac{\partial}{\partial y_{nn}} & \frac{\partial}{\partial x_{nn}} \\
\end{array} \right). ~~~
\eea
Clearly they are decomplexifications of matrices $Z$ and $D$ defined by
formula \ref{ZFnot1} above. (We assume that   $Z$ is defined as: $Z_{ab} = x_{ab}+iy_{ab}$
and $D$ is defined as: $D=\partial_{Z_{ab}} = 1/2(\partial_{x_{ab}}- i \partial_{y_{ab}})$).

Our first new result is: 

{\bf Theorem \pref{ThCap} }  (section \ref{SSSectCap}).  

{\em 
\bea \label{CapInt1}
{\det}^{\mathrm{col}} ( Z^{\RR} (\DD^{\RR})^t + {CorrTriDiag})= {\det}
 ( Z^{\RR}) {\det} ((\DD^{\RR})^t),
\eea
where there are two options for the matrix $CorrTriDiag$ , 
the first one is:
\bea \label{fmlCorrTriDiag}
CorrTriDiag  = \left( \begin{array}{ccccccccc}
n-1 +1/4 & i/4 & 0 & 0 & 0 & 0 & ... & 0 & 0 \\
i/4 & n-1 - 1/4 & 0 & 0 & 0 & 0 & ... & 0 & 0 \\
 0 & 0 & n-2 +1/4 & i/4 & 0 & 0 & ... & 0 & 0 \\
 0 & 0 & i/4 & n-2 - 1/4 &0 & 0 & ... & 0 & 0 \\
 ... & ... & ... & ... & ...  & ... & ... & ... & ... \\
 0 & 0 & 0 & 0 & 0 & 0 & ... & 1/4 & i/4 \\
 0 & 0 & 0 & 0 &0 & 0 & ... & i/4 & -1/4 
\end{array} \right). ~~~
\eea
The second option for $CorrTriDiag$ is the same formula with substitution $i$ to $-i$,
i.e. complex conjugate to the first one.
}

{\bf Theorem \pref{ThASymCap}}  (section \ref{SSSectASymCap}).

{\em    
 ("Decomplexification" of Turnbul symmetric and  Howe-Umeda Kostant-Sahi ant-symmetric identities).
The same formula \ref{CapInt1} holds if $Z^\RR$ and $\DD^\RR$ are "decomplexifications"
of either symmetric or anti-symmetric matrices made of $z_{ij}$ and $\partial_{z_{ij}}$.}

So, the correcting term, which was 
$\mathrm{diag}(n-1,n-2,...,1,0)$ in the classical identity changed by the rule that each number $k$ is substituted by the two by two matrix:
\bea
k ->  \left( \begin{array}{ccccccccc}
k +1/4 & i/4  \\
i/4 & k - 1/4
\end{array} \right). ~~~
\eea
This is a result of some calculation without clear explanation, more precisely it is natural to expect that 
correcting term has the form $[k +a ~,~ b ; ~~ c ~,~ k +d]$ with determinant equal to $k^2$, but this requirement does not fix $a,b,c,d$.  The fact
that correction is not so simple - not like $diag(k,k)$ seems surprising for us.
Let us mention that the matrix of the form 
$\left( \begin{array}{ccccccccc}
1/4 & i/4  \\
i/4 & - 1/4
\end{array} \right)$ is matrix of eigenvectors of any decomplexified 2x2 matrix, may be this
can give key to explanation, but for the moment it is unclear. See also remark \pref{EigRem}.

Our next result is an extension of the decomplexified Capelli identities to the case of rectangular matrices,
i.e. non-commutative analog of the Cauchy-Binet formula. 

{\bf Theorem \pref{ThRectCapTur} } (section \ref{SSSectRectCapTur}).  

{\em  
Let $Z,D$ be square $n\times n$ matrices which are
  
either

(Capelli's case) 
$Z$ is a matrix with $z_{ij}$ at place  $(ij)$, and $\DD$ with elements $\partial_{z_{ij}}$ respectively
(see formula \ref{ZFnot1} above)

or 

(Turnbull's case) $Z,D$ are  
symmetric matrices  defined as: 
 $z_{ij}$ $i\le j$ on both positions $(ij)$ and $(ji)$ and 
$\DD$ is the symmetric matrix with  $\partial_{z_{ij}}$ out of the diagonal and $2 \partial_{z_{ii}}$ on the diagonal

For any multi-indexes $I=(i_1,...,i_{r})$,
$J=(j_1,...,j_{r})$ the following identity is true:
\bea
{\det}_{2r\times 2r}^{\mathrm{col}} ( (Z ~\DD^t)^{\RR}_{IJ} +Q^{\RR} {CorrTriDiag})= 
\sum_{L=(l_1<l_2<l_3<...<l_{2r})}
{\det}_{2r\times 2r} ( Z^{\RR}_{IL})
 {\det}_{2r\times 2r} ((\DD^{\RR})^t_{LJ}).
\eea
}

Where $CorrTriDiag$ is defined by  formula
\ref{fmlCorrTriDiag} above.
And matrix $Q$ is  defined:
\bea
(Q)_{ab} = \delta_{i_a j_b} .
\eea

Apparently, a similar result is true in the antisymmetric case. We formulate conditional theorem about the
antisymmetric case in section \ref{SSSectRectAntiSym}. 

Next we propose certain "decomplexification" for the Capelli identities
discovered in \cite{CSS08} (some generalization and alternative proofs were given
in \cite{CFR09}). 

{\bf Theorem \pref{CCSCap} } (section \ref{CCSCapSSSect}).  

{ \em
Assume  $M,Y,Q$ are square $n\times n$ matrices, such that: \\
either\\
(1) $M,Y,Q$ satisfy CSS-condition $ \forall l,j,r,p:~~ [Y_{lj}, M_{rp}] = \delta_{lp} Q_{rj}$ 
 \\
(2) $M$ is a Manin matrix (see e.g. subsubsection \pref{SSSectManinMatrReminder} )\\
or\\
(1') $M,Y$ satisfy TCSS-conditions $\forall i,j,k,l:~~ [M_{ij}, Y_{kl}] = -h(\delta_{jk}\delta_{il} + \delta_{ik}\delta_{jl})$ and $h$ is central, matrix $Q$ is defined as $Q=h~Id$ \\
(2') $M$ is a symmetric matrix with commuting entries\\
and in both cases CSS and TCSS we require:\\
(3) elements in the same column of $Y$ commute among themselves i.e. $\forall a,b,k: [Y_{ak}, Y_{bk}]=0$; \\
(4) For any two matrix elements $\alpha, \beta$ of matrices $M,Y,Q$, it is true that $[\alpha, \bar \beta]=0$, i.e.
all elements commute with complex conjugated elements;  \\ 
{\bf then}
\bea
{\det}^{\mathrm{col}}_{2n\times 2n} ( M^{\RR} Y^{\RR} + Q^\RR\mathrm{CorrTriDiag})= {\det}^{\mathrm{col}}_{2n\times 2n} ( M^{\RR}) {\det}^{\mathrm{col}}_{2n\times 2n} (Y^{\RR}).
\eea
}
Where $CorrTriDiag$ is defined by the formula
\ref{fmlCorrTriDiag} above.

\subsubsection{Holomorphic factorization results}

Actually, the theorems above follow quite easily  from 
known Capelli identities and the  
following result, which  is an analog of the well-known fact that 
$\det(M^\RR) = \det(M) \det(\bar M)$ for the determinants
of the "decomplexified" commutative matrices.

{\bf Theorem \pref{HolFactCapel}} {\em  (section \ref{SSectHolFactCapel} ).  

 Consider elements $E_{ij} = E_{ij}^{Re} + i  E_{ij}^{Im}$, satisfying commutation relations:
$[ E_{ij} , E_{kl} ] = E_{il} \delta_{jk} - E_{kj}\delta_{li} $ 
i.e. the same relations as for matrices with zeros everywhere except $1$ at the position $(ij)$.  Moreover impose relation $[ E_{ij}, \bar  E_{kl}  ]=0$, i.e. 
all elements commute with complex conjugated elements.
 Consider the $2n\times 2n$ matrix $E^{\RR}$: 
\bea
E^{\RR} = \left( \begin{array}{ccccc}
E_{11}^{Re} & E_{11}^{Im} & ... & E_{1n}^{Re} & -E_{1n}^{Im} \\
-E_{11}^{Im} &E_{11}^{Re} & ... & -E_{1n}^{Im} & E_{1n}^{Re} \\
... & ... & ... \\
E_{n1}^{Re} & E_{n1}^{Im} & ... & E_{nn}^{Re} & -E_{nn}^{Im} \\
-E_{n1}^{Im} & E_{n1}^{Re} & ... & -E_{nn}^{Im} & E_{nn}^{Re} \\
\end{array} \right), ~~~
\eea

then
\bea
{\det}^{\mathrm{col}}_{2n\times 2n} (E^{\RR} + CorrTriDiag ) = \\ = 
{\det}^{\mathrm{col}}_{n\times n} (E +  diag(n-1, n-2, ...,1,0) ) 
{\det}^{\mathrm{col}}_{n\times n} (\bar {E}  + diag(n-1, n-2, ...,1,0) ), 
\eea
here $E$ is an $n$ by $n$ matrix with $E_{ij}$ on the position $(ij)$,
respectively $\bar {E}$ with $\bar E_{ij} = E_{ij}^{Re} - i  E_{ij}^{Im}$.
(Clearly $E^\RR$ is  decomplexification of $E$).
The matrix $CorrTriDiag$ is defined by the formula \ref{fmlCorrTriDiag} above.
}

Actually, it is more convenient to prove a more general result given below.
It allows us to see more clearly what properties are used at each step.
It is the main theorem of the paper.

{\bf Theorem \pref{MainTheor} }(section \ref{MainThSect}).
 {\em  

Consider $n\times n$ matrices $C,Q$,
for any two matrix elements $\alpha, \beta$ of matrices $C,Q$, it is true that $[\alpha, \bar \beta]=0$, i.e.
all elements commute with complex conjugated elements.
Assume also that their elements satisfy the following relations
(the relations are column-wise i.e. NO relation between elements from different columns):
\bea
\forall i,j,k: ~~[C_{ik}, C_{jk}] = C_{ik} Q_{jk} -  C_{jk}Q_{ik},  
~~~ [C_{ik}, Q_{jk}] = [C_{jk}, Q_{ik}],
~~~ [Q_{ik}, Q_{jk}] = 0, 
\eea
{\bf then} the following analog of holomorphic factorization holds:

for any central elements $d_i$: 
\bea
{\det}^{\mathrm{col}}_{2n\times 2n} (C^{\RR} + Q^{\RR} CorrTriDiag ) = \\ =
{\det}^{\mathrm{col}}_{n\times n} (C + Q diag(d_n, d_{n-1}, ...,d_{2}, d_{1} ) ) 
{\det}^{\mathrm{col}}_{n\times n} (\bar C +\bar Q diag( d_n,  d_{n-1}, ...,  d_{2}, d_{1} ) ) ,
\eea
where there are two options for the matrix $CorrTriDiag$ , 
the first one is:
\bea
CorrTriDiag  = \left( \begin{array}{ccccccccc}
d_n +1/4 & i/4 & 0 & 0 & 0 & 0 & ... & 0 & 0 \\
i/4 & d_n - 1/4 & 0 & 0 & 0 & 0 & ... & 0 & 0 \\
 0 & 0 & d_{n-1} +1/4 & i/4 & 0 & 0 & ... & 0 & 0 \\
 0 & 0 & i/4 & d_{n-1} - 1/4 &0 & 0 & ... & 0 & 0 \\
 ... & ... & ... & ... & ...  & ... & ... & ... & ... \\
 0 & 0 & 0 & 0 & 0 & 0 & ... & d_{1}+1/4 & i/4 \\
 0 & 0 & 0 & 0 &0 & 0 & ... & i/4 & d_1 -1/4 
\end{array} \right). ~~~
\eea
The second option for $CorrTriDiag$ is the same formula with substitution $i$ to $-i$,
i.e. complex conjugate to the first one. 
}

 We also discuss Cayley identities, but
these results are fairly straight-forward, or open questions, so we will not mention them
in introduction, just referring to the main text.

\subsection{ The Capelli identities and around}
Here we  describe the relation of the Capelli identities with other questions and 
give some brief overview of the related works. 
We have 4 subsections 1) describe relation with the Wick ordering 
2) describe applications
in representation theory 
3)
present observation that Capelli determinant equals to the Duflo image of the commutative determinant,
thus providing  conceptual meaning for the Capelli determinant
4) give brief overview of the Capelli identities.

\subsubsection{Capelli means: Wick ordered quantization is "great". Open issues\label{WickIntrSect}}

Let us give another point of view on Capelli identities 
quite
inspiring in certain questions and relates them with the general problem of ordering in
quantization theory.
Commutative algebras are related to classical physics, and non-commutative algebras
to quantum physics. Hamiltonians of classical systems are typically some
functions of commuting variables $q_i, p_i$.
While Hamiltonians of quantum systems are functions of non-commuting variables $\hat q_i, \hat p_i$,
such that $[\hat q_i, \hat p_j] = \delta_{ij}$ (we omit Planck constant)
and $[\hat q_i, \hat q_j] = [\hat p_i, \hat p_j] = 0$.
To get quantum Hamiltonian from classical Hamiltonian, one needs 
somehow to resolve the "ordering problem" (i.e. if classical expression contains
$qp$ it may correspond to $\hat q \hat p$ or $\hat p \hat q$ or $1/2(\hat p \hat q +\hat q \hat p)$ etc).
The simplest recipe is "Wick ordering" i.e.  put the $\hat q_i$ lefter than
$\hat p_i$:
\bea
Wick: ~ \CC[q_i, p_i]\to \CC[\hat q_i,\hat p_i],\\
Wick(\prod q_i^{k_i} p_i^{n_i}) = (\prod \hat q_i^{k_i} \hat p_i^{n_i}) .
\eea

The choice of ordering is not unique.
And the important question is: if we have some good property for the classical
system (e.g. system is integrable) will it be preserved on the quantum level
if we use certain ordering ? It is believed that 
no ordering will work
well for all situations. 
The correct construction of a quantum system from a classical one
remains somehow an art.
In this respect, the claim below, that in certain situations Wick ordering (which is very simple) works very well, comes as a surprise.

Now we are ready to formulate the first message of this subsection.
{\bf For the classical systems whose Hamiltonians are related with 
Capelli identities (e.g. Gaudin model)
the Wick ordering is expected to work very well}.
For example for the Gaudin model the identity \cite{MTV-Cap}
implies that Wick ordering of classical Hamiltonians produce quantum
commuting Hamiltonians, which is a highly non-trivial fact, and the only existing proof
uses the Capelli identity and Talalaev's theorem \cite{T04}.

Let us reformulate the Capelli identities
with the help of the Wick ordering.
Let $Z,D$ be the same as above - formula \ref{ZFnot1}.
In our paper we do not use $q,p$ notations for variables, but $z_{ij}$ plays
role of $q$, variables $p_{ij}$ play role of $p$, respectively
$\partial_{z_{ij}}$ play role $\hat p$ and $\hat z_{ij} =z_{ij}$.
So Wick ordering in our notations looks as - put $z_{ij}$ to the left of 
 $\partial_{z_{ij}}$:
\bea
Wick(\prod z_{ij}^{k_{ij}} p_{ij} ^{n_{ij}}) = 
\prod z_{ij}^{k_{ij}} \partial_{z_{ij}} ^{n_{ij}}.
\eea 

Formulations of the Capelli with Wick:

\bea
{\det}^{\mathrm{col}}(ZD^t+diag(n-1,n-2,...,1,0) ) = Wick(det(ZP^t)) .
\eea
The standard Capelli follows immediately: $Wick(det(ZP^t))=
Wick(det(Z)det(P^t))=Wick(det(Z))Wick(det(P^t))=det(Z) det(D^t)$.

\bea
perm^{col} (Z^tD - diag(n-1,n-2,...,1,0) ) = Wick(perm(Z^tP)) .
\eea
Again, it follows from the known Capelli identities for permanents
\cite{Wil81},\cite{Na91}, \cite{CFR09} in the same way as above. Note that
in commutative case there is Cauchy-Binet formula for permanents (e.g. proposition 13 page 49 \cite{CFR09}).

\bea
Immanant_{quantum} (ZD^t ) = Wick(Immanant_{standard} (ZP^t)) .
\eea
Again it follows from the known Capelli identities for immanants\footnote{
Classical immanants are parametrized by representations $\pi$ of $GL_n$ and 
$Immanant(g) = Tr_{\pi} (\pi(g))$ for matrix $g \in GL_n$.}
\cite{Wil81}, \cite{Ok96-1}. Construction of "quantum immanants"
in \cite{Ok96-1} is also called "fusion procedure" in physics literature (e.g.
\cite{KSZ07}
\footnote{There are quite a lot physicist's activities on quantum immanants (e.g. "fusion") which is in somewhat "parallel universe" to mathematician's activities.
  Immanants for diagonal matrices are precisely
Schur functions which satisfy plenty relations e.g. Giambelli, Jacobi-Trudi.
There are quite many physicist's papers deriving "quantum" analogues of these relations, which are sometimes called "Cherednik-Bazhanov-Reshetikhin" formulas
see e.g. \cite{AKLTZ11} and references there in.
These considerations are important in AdS/CFT correpondence - hence
 they became quite a hot topic in physics literature (see \cite{AdSCFTInteg}
 chapter III.7).
}).

The \cite{MTV-Cap} identity (see also \cite{CF07}) is even more striking since
there is no need for any correction terms on the left hand side.
\footnote{here we give a simplified example of \cite{MTV-Cap},\cite{CF07} result
which is actually equivalent to the standard one. More general
version is same statement but for more general matrix:
 $L(t) = \partial_t - A - Z(t-B)^{-1} D^t$, where $A,B$ are constant $\CC$-valued 
 matrices and $X,D$ can be rectangular} 
One needs to introduce auxiliary  commuting variables $t$ and $\lambda$
with the agreement that $Wick(\lambda) = \partial_t$. 
\bea \label{MTVCAPintro}
{\det}^{\mathrm{col}}(\partial_t - ZD^t / t) = Wick(det( \lambda - ZP^t / t)).
\eea
So Wick ordering gives us not only commuting Hamiltonians, but also important differential operator in $t$ (called "GL-oper" in Langlands theory, 
"Baxter's T-Q relation " or "Sklyanin's quantum spectral curve" in
integrability (\cite{CT06-1})). 
Let us remark that $\partial_t - ZD^t / t$ is an example of Manin matrix,
and the idea to put $\partial_t$ is a particular case of  Talalaev's \cite{T04} idea.

So the second message of  this subsection is that  all the identities
above have the form
\bea
F_{quantum} (Wick(ZP^t)) = Wick ( F_{classical} (ZP^t)),
\eea
where $F_{classical}$ is  some invariant function (determinant, permanent, etc...)
and $F_{quantum}$ is its appropriately understood quantum version i.e.
some correcting terms introduced or fixed row-column order, etc. 
And the statement has the form that {\bf "Wick and $F$ commute" } modulo
these corrections.

For certain reasons we believe that the left hand side quantum expressions
(e.g. Capelli determinants, permanents, immananats...) are "correct quantizations"
of the corresponding classical expressions (because they are very good from representation theory and integrability points of view - they give central elements,
they preserve integrability and so on). 
The moral of these identities is: {\bf these "correct quantizations" 
are just Wick orderings of the corresponding classical expressions}.

We think that the situation above is specific to Capelli's setup,
but it  might be that it is not the only example.
Understanding when Wick ordering  is a "correct quantization"
seems quite an interesting issue.

{\bf Open questions.}

The open questions formulated in \cite{CF07}: whether the  identity \ref{MTVCAPintro} holds true if we substitute $det$ with permanent and trace.

The general questions - can the identities above be extended to the case
of symmetric and antisymmetric matrices $Z$ and $D$ ?
We may also consider their extension to rectangular matrices.

About the first one  we know from Turnbull (for symmetric case) and Howe-Umeda Kostant-Sahi 
(for antisymmetric n-even, (the n-odd is an open conjecture in \cite{CSS08}
page 36 bottom conjectures 5.1, 5.2)) that the identity  is true.
The second identity is actually also known for antisymmetric matrices
from Turnbull (see \cite{FZ93} last theorem page 10 (without proof)  and \cite{CSS08} proposition 1.5 page 9, and it also follows combining  \cite{CFR09} theorem 13 page 50 and proposition 13 page 49).

The other cases remain open to the best of our knowledge.

And finally one can try to do "decomplexification" and "dequaternionifications" of the identities above, similar to what is done in the present paper.

Wick ordering has been generalized to quantization of symplectic manifolds \cite{DLS01},
it would be interesting to explore orbits and their quantization and try to get generalize Capelli identities
in this set up.

\subsubsection{Why matrices with non-commutative entries are interesting \label{sectWhy} }

{\bf Reason 1.}
Matrices with noncommutative elements appear in representation theory and quantum integrable systems theory.
Roughly speaking the reason is the following: one wishes to consider certain  operators  (e.g. commuting Hamiltonians in integrability theory or elements of centers,
intertwiners, etc.   in representation theory) or, more generally, elements of some non-commutative algebras. It turns out that explicitly written formulas for them are rather complicated, but {\bf the insight is that complicated operators can be constructed as determinants, traces, etc... of some matrices with noncommutative entries, while the matrices themselves have much simpler form and more easier to deal with than the operators of interest).}

{\bf Reason 2.}
The second reason for us is actually an illusive dream. The non-commutative elements  can be matrices themselves. So we get usual matrices, with additional
structure of partition to blocks. If these blocks satisfy certain specific commutation relations, then  many useful facts of commutative linear algebra can be extended to them
and $det(A)det(B) = det(AB+correction)$ is just an example.
In particular, 
{\bf if such block-matrices appear in some practical application, then standard algorithms (matrix inversion, svd, etc.) can be improved  
to work faster.} The problem is that now it seems hard to imagine real life application where such block-matrices may appear.

{\bf Examples for reason 1.}


Probably the main  example is the quantum group theory, which originated as
a mathematical framework for the Yang-Baxter equation and related $RTT=TTR$
matrix equations. So non-commutative matrices are intimately related to quantum
groups. But let us provide examples from the more classical representation theory,
which are quite related to the present paper. Let us also refer to  survey by A. Molev \cite{Molev02} and references therein, and recent work \cite{AG11} for other examples.

{\bf Proposition (Capelli).}
{\it 
Consider 
\bea \label{CapCenterFml}
{\det}^{\mathrm{col}}_{n\times n} (E +  diag(n-1, n-2, ...,1,0) +t ) 
= \sum_k C_k t^k.
\eea
All elements $C_k$ belong to the center of $U(gl_n)$. (Actually they freely generate the center.) Here matrix $E$ is defined as usually as $E_{ij}$ stands on position $(ij)$, where $E_{ij} \in gl_n\subset U(gl_n)$ is the matrix with zeros everywhere except $1$ at the position $(ij)$. 
}

The easiest proof is based on the Capelli identity: consider for simplicity only the determinant $C_0$. Capelli identity states $C_0 = det(Z)det(D)$. To check that the element is central is the same as to check
that it is invariant with respect to action of $GL_n$ by conjugation, this is easy to see:
$det(g Z g^{-1}) det(  g Z g^{-1} ) = det(g) det(Z) det(g)^{-1} det(  g) det( Z) det(g)^{-1} =
det(Z)det(D)$. To get similar proof for $C_i$ we should consider Capelli identities for rectangular matrices,
it is a little longer, so let us omit it.

The second example is related to 
Talalaev's formula and the center of $U_{crit}(gl_n[t,t^{-1}])$.

Consider Lie algebra $gl_n[t,t^{-1}]$ and its universal enveloping  $U_{crit}(gl_n[t,t^{-1}])$
(actually several technical details are necessary - central extension, putting central element to be 
equal to some specific number "critical level" and making the local completion).
Denote as $L(z)$ the matrix:

\bea \label{L-full}
L(z)=\sum_{i=-\infty ... \infty }
\frac{1}{z^{i+1}}
\left(\begin{array}{ccc}
E_{1,1}t^i & ... & E_{1,n}t^i \\
... &  ... & ... \\
E_{n,1}t^i & ... & E_{n,n}t^i
\end{array}\right) 
\eea
Here, $z$ is a formal parameter, $E_{ij}t^k$ are elements of $gl_n[t,t^{-1}]$ (more precisely  elements of  $U_{crit}(gl_n[t,t^{-1}])$).
The following theorem is based on D. Talalaev's breakthrough result \cite{T04}.

{\bf Theorem \cite{CT06-1}, \cite{CM} }  
{\em The elements $C_{i,k} \in U_{crit}(gl_n[t,t^{-1}])$ defined by Talalaev's type formula below 
freely generate the center of $U_{crit}(gl_n[t,t^{-1}])$
\bea
:{\det}^{\mathrm{col}}(\partial_z -L(z)): = \sum_{i=0...n,k=-\infty ... \infty} \partial_z^i z^k  C_{i,k},
\eea
here $:...:$ denotes the standard "normal ordering" ($:a(z)b(z):= a(z)_{+}b(z) + b(z)a(z)_-;$,
where $a(z)_+=\sum_{i \ge 0} a_{-i-1} z^i;  a(z)_-=\sum_{i < 0} a_{-i-1} z^i$).
}

The result has been extended to the super case $gl(n|m)[t,t^{-1}]$ \cite{MR09},
and to  other classical Lie algebras ($B,C,D$-series) in \cite{M11}.
Related results are in \cite{RST} for elliptic algebras
and q-analogs \cite{CFRS}.
Let us mention the paper \cite{MTV05-12}, where Talalaev's formula is used to settle some long-standing conjecture in  real algebraic geometry.

Let us also mention that the matrix $(\partial_z -L(z))$ is a Manin matrix
with respect to the normal ordered product (\cite{CF07}, \cite{CM}). 

And finally remark that if we restrict the Lax matrix from $gl_n [t,t^{-1}]$
to $gl_n$ itself, i.e. just  preserve only the term corresponding to $i=0$ in
 formula \ref{L-full}, then Talalaev's formula will give elements of the center of $U(gl_n)$,
and they will coincide with the Capelli elements (\cite{MTV-Cap}).
In that sense Talalaev's formula provides correct generalization to $gl_n[t,t^{-1}]$ of the Capelli elements for $gl_n$.

\subsubsection{Capelli determinant = Duflo (determinant) \label{CapDufloIntrSect}}

The observation below seems to be new. 
Briefly speaking we claim that the Capelli determinant equals to the image under the Duflo map for the Lie algebra $gl_n$ of the standard determinant of matrix with commuting entries. 
Duflo map \cite{Duf77}   is something amazing and general, so the claim provides a conceptual meaning of the Capelli determinant which has been introduced somewhat ad hoc
in nineteenth century. Let us go on with details. We refer to \cite{Kir99}, \cite{CR11} for information on Duflo map.

Consider square $n \times n$ matrix $E$ with elements $E_{ij}$.
Assume elements $E_{ij}$ satisfy commutation relations: 
\bea \label{glnComRel}
[E_{ij}, E_{kl}] = \delta_{jk} E_{il} - \delta_{il} E_{kj},
\eea
i.e. the same relations which are satisfied by the matrices which have zeros everywhere except 1 at position $(ij)$ - i.e.  basic generators of Lie algebra $gl(n)$.
Under this assumption, as is well-known \footnote{ See formula \pref{CapCenterFml} and around}, 
for all $u\in \CC$ Capelli determinants ${\det}^{\mathrm{col}}(E+diag(n-1,n-2,...1,0)+u~Id)$ commute will all $E_{ij}$ i.e. they are central element in algebra generated by $E_{ij}$.

Consider universal enveloping algebra $U(gl_n)$, i.e. algebra  generated by $E_{ij}$ with the only relations \ref{glnComRel}. Consider also symmetric algebra $S(gl_n)$ which is commutative polynomial algebra $\CC[E_{ij}]$. For any finite-dimensional Lie algebra $g$ over field of characteristic zero
M. Duflo constructed the map: 
\bea
Duflo: S(g) \to U(g), 
\eea
which is being restricted on $S(g)^{g}$ gives {\bf isomoprhism of commutative algebras
$S(g)^{g}$ and  $ZU(g) = U(g)^{g}$ }. Also Duflo map is\footnote{on the whole $S(g)$ without restriction to
$S(g)^{g}$ }   isomorphism of graded vector spaces and $g$-modules. Moreover it is  identity map on associated graded algebras $Id = gr \circ Duflo: S(g) \to Gr(U(g)) = S(g)$. 
Let us emphasize that non-trivial part of the statement is isomorphism of commutative subalgebras.

{\Prop Consider matrix $E$ with elements $E_{ij}$. Slightly abusing notations, we will look at $E_{ij}$ as generators of two algebras: first algebra -  commutative $S(gl_n)$ and second algebra - non-commutative $U(gl_n)$.  Consider commutative algebra $S(gl_n)=\CC[E_{ij}]$ and calculate $\det(E)$ as standard determinant of matrix with commutative elements.
Consider non-commutative algebra -  universal enveloping $U(gl_n)$ and calculate Capelli determinant
${\det}^{\mathrm{col}}(E+diag(n-1,n-2,...1,0)- (n-1)/2Id)$ taking into account relations \ref{glnComRel}.
Then:    
\bea 
Duflo( \det( E) )  = {\det}^{\mathrm{col}}(E + diag(n-1,n-2,...1,0) - (n-1)/2Id) ,
\eea
i.e. the image of the standard determinant under the  Duflo map equals to the Capelli determinant.
\footnote{Pay attention that in Capelli identity we have 
${\det}^{\mathrm{col}}(E+diag(n-1,n-2,...1,0))$, while here: ${\det}^{\mathrm{col}}(E+diag(n-1,n-2,...1,0)- (n-1)/2Id)$, the difference is $- (n-1)/2Id$.}
}

{\Rem ~}   $diag(n-1,n-2,...1,0) - (n-1)/2Id =  diag( (n-1)/2, (n-3)/2, ...-(n-3)/2,-(n-1)/2))$.
We get trace zero diagonal correction, which is rather natural.

\PRF\footnote{Our motivation to guess this fact was rather indirect and related to \cite{CR04}:
from this conjecture (which was later proved in \cite{Cat07})  we can conclude that Duflo determinant should vanish in representations
corresponding to orbits of dimensions less than generic orbit - roughly speaking standard determinant 
vanishes on matrices of rank less than full and Duflo map respects this - the same can be directly seen for the Capelli determinant,
moreover their symbols coincide - hence the difference can be only the minors of lower degree - but they do not vanish on smaller dimensional orbits, hence the difference is zero.}  Observe that both left hand side and right hand side belong to the center of $U(gl_n)$.
Left hand side - due to Duflo, right hand side essentially due to Capelli.
It is known, that to check that two elements of $ZU(gl_n)$ coincide it is enough to check that their values
in irreducible highest weight $(\lambda_1,...\lambda_n)$ representations coincide,
i.e. their Harish-Chandra images coincide.  
It is rather easy to see that image of the determinant on Capelli side equals to 
$ (\lambda_1 + (n-1)/2)  (\lambda_2 + (n-3)/2) ... (\lambda_n - (n-1)/2)$, because we can calculate the
desired value on the highest weight vector itself, and so all upper diagonal elements of $E$ act by zero on it and 
so determinant equals to the product of diagonal terms, which equals to expression above on the highest weight vector.
On the other hand by Duflo it is known how to calculate Harish-Chandra image of any element obtained by the Duflo map. Consider any element $A$ in $S(gl_n)^{gl_n}$ it can be identified with some symmetric polynom $p(\lambda_i)$, then $Harish-Chandra(Duflo(A))=p(\lambda_1 + (n-1)/2, \lambda_2 + (n-3)/2, ... , \lambda_n - (n-1)/2   )  $.
Applying Duflo's result to $\det(E)$ which corresponds to symmetric function $\prod_i \lambda_i$
we get the same result as for Capelli determinant, so proposition is proved.
\BX

In the same way we can prove similar result for Capelli permanents, and then MacMahon-Wronski
relation between Capelli  permanents and Capelli determinants follows
 from the same relation in commutative case - just because Duflo map is homomorphism.
 This gives alternative proof of the results \cite{Umeda03}.
It would be interesting to clarify the relation between quantum immanants \cite{Ok96-1}
and Duflo images of classical immanants.
In particular in \cite{BB06} certain Capelli type determinants has been constructed,
it might be they can also be equal to the Duflo images of corresponding classical expressions.


M. Kontsevich \cite{Kon97} generalized the Duflo map to quantization of an arbitrary Poisson manifold,
this stimulated much interest in Duflo map - see \cite{CR11}.
He showed that Poisson center of Poisson algebra $Fun(M)$ is isomorphic to the center of 
the non-commutative algebra $Quantization(Fun(M))$ for any smooth manifold $M$.
So if some element $d$ belongs to the Poisson center of  $Fun(M)$
we should consider $Duflo-Konstsevich (d)$ as "correct quantization of an element $d$".
Respectively proposition above shows that Capelli determinant is
correct quantization of a standard determinant.
Highly intriguing is to consider the case of loop algebras, try to define Duflo map for them (standard construction does not work, since algebra is infinite-dimensional) and try to relate 
with Talalaev's formula.

To conclude Capelli elements in the centers are down-to-earth, simple-to-work, but somewhat ad hoc,
elements, on the other hand Duflo map gives general and conceptual way to construct such elements,
but it is difficult to write down them explicitly in some PBW-basis. The equality between the two elements
gives a pleasant conclusion that we have elements which are both easy-to-work and conceptually clear.
It would be interesting to generalize this result for other Lie algebras and not only for  determinants,
permanents.

\subsubsection{Brief overview of the Capelli identities\label{OverCapSect}}

Let us briefly mention surrounding works on the Capelli identities, which are known to us
to some extent. This brief survey is not complete.
We refer to \cite{Um98} and to Wikipedia article
\cite{WikiCapelli} for further information. 
Current version of Wikipedia article
has been mainly
written by the present author. Wikipedia is a collaborative resource,
so let me express a hope that some one will participate in improving the article.

The "pre-history" period -  let us mention only Alfredo Capelli himself,
Herbert W. Turnbull, the book by Herman Weyl "The Classical Groups, their Invariants and Representations", G.C. Rota  \cite{GCRota69} as referred in \cite{Wil81},\cite{CW81}.

The modern interest
strongly motivated by the works by Roger Howe who related the identities with the theory
of dual pairs \cite{H89}.
Several papers appeared at the beginning of 90-ies 
\cite{HU91}, \cite{KS91}, \cite{FZ93}, etc...
Important bunch of works by the
Russian team (members of  A.A. Kirillov's seminar) A.Molev, M. Nazarov, A. Okounkov and G. Olshanski  appeared in the middle of 90-ies.
These developments were partly motivated by the discovery that
technique originating in integrable system theory and quantum groups (Yangians)
can be applied to the subject, in particular Capelli determinant can be obtained as Yangian's quantum determinant. There were obtained generalizations to immanants ("higher Capelli identities"), generalizations to semi-simple Lie algebras and to their super-analogues - we refer to surveys  \cite{MNO}, \cite{Molev02}, \cite{WikiCapelli} and references therein. Japanese researchers (T. Umeda, his colleagues and students) extended the Capelli identity
to quantum groups (see references in \cite{Um98}, \cite{WikiCapelli}). 
Later on they concentrated 
the classical Lie algebras.
Many dual pairs has been considered, results extended from determinants
to permanents, traces and immanants, various non-commutative linear algebra
theorems has been discovered (Cayley-Hamilton, Cramer theorems, Newton
formulas, etc.);
simple proofs for some classical identities has been worked out.
(See papers by T. Umeda, T. Hashimoto, M. Itoh, K. Kinoshita, K. Nishiyama,  A. Wachi, M. Wakayama  referred in \cite{Um98}, \cite{WikiCapelli}).
The short and very readable paper by Toru Umeda \cite{Um98AMS} 
is highly recommended.

Probably around 2000 one might have feeling that the story may  be not ended, but the landscape is quite seen. The generalizations come in several directions:
1) 
various dual pairs for different super-q-semisimple Lie algebras;
2) determinants
can be changed by permanents, pfaffians, traces, immanants, etc., 
3) not only the identity $det(AB)=det(A)det(B)$, but 
other facts of linear algebra (Cayley-Hamilton, Cramer, Newton, ...)
become of interest for  Capelli type  matrices with non-commuting elements.
(See \cite{CFR09} for extensive list of references on "non-commutative linear algebra").


But after 2006 there appeared several unexpected and remarkable achievements in the area.
In \cite{MTV-Cap} remarkable generalization of the 
Capelli identity related to  integrability (see also \cite{CF07})
has been found. It should be thought of as a kind of extension of the
Capelli theory from finite-dimensional Lie algebras to their loop algebras.
So the questions how to generalize previously known results to loop algebras,
understanding their relation to integrability and Langlands correspondence
arise.
In \cite{CSS08} unexpected generalization appeared - previously all identities
has been related to semisimple Lie algebras or their q-super-deformations - but the
identity found in \cite{CSS08} has nothing to do with this.
It has purely algebraic form claiming that if certain commutation
relations between matrix elements holds true then there exists
an analog of the formula $det(A)det(B)=det(AB)$, but in the right
hand side appeared Capelli-like correction $det(AB+Q)$,
either $A$ or $B$ is \MM ~in these formulas.
Some generalizations and alternative proofs were proposed in \cite{CFR09}.
Finally in 2011 another unexpected formula has been proposed
in the paper \cite{AH11}, the present paper is  motivated by it.
The formula claims Capelli type identity for $4\times 4$ matrices
$det(H^{\RR}) det(D^{\RR}) = det( H^{\RR} D^{\RR}+Correction)$.
The matrices $H^{\RR}, D^{\RR}$ have a form of "dequaternionification" i.e. natural representation
of quaternions by real valued $4\times 4$ matrices,
the matrix $H^{\RR}$ contains $x_1, x_2,x_3,x_4$, while $ D^{\RR}$ contains
$\partial_1, \partial_2,\partial_3,\partial_4$.
So in the present paper we consider similar situation but for "decomplexification"
and our work treats the general case of $n\times n$ matrices.
Also we get various generalizations like rectangular, (anti)-symmetric, etc.
identities.

Last, but not least let us mention brilliant, but undeservedly forgotten
paper from the  combinatorial community - paper by S.G. Williamson
\cite{Wil81}
who discovered generalization of the Capelli identities to permanents and immanants,
somehow foreseeing the works by A. Okounkov and M. Nazarov decade later.
His approach  avoids "fusion procedure" which was the key ingredient of
the later works and constructs "higher Capelli" identities purely by combinatorial
means. Let us also mention the paper \cite{CW81} where 
remarkable generalizations of the Cayley type identity has been found.

\subsection{Comments on results, paper organization, proofs etc.}
We have 4 subsections here: 1) explain what we prove first, what we prove second and why
 2) describe organization of the paper 3) explain the main idea of the proof 
 4) make further rather informal comments.

\subsubsection{Brief summary of results and their inter-relations}
We have two kinds of main results - first are about Capelli identity
($det(A)det(B) = det(AB+correction)$),
another about holomorphic factorization ($det(A^\RR+cor)= det(A+cor)det(\bar A+\widebar{cor})$). We also discuss Cayley identities, but
these just for completeness - these results are fairly straight-forward or open.
All Capelli type results follows from holomorphic factorization results and known 
previously Capelli identities.
So let us first discuss factorization results.
Our key general result - is theorem 
 \ref{MainTheor} in section \ref{MainThSect}. It gives certain general
conditions when factorization result holds true.
Next we easily check that Capelli type matrix  satisfies these conditions,
hence get factorization result for it. It is necessary for Capelli identities.  

Capelli results are splitted into three - for square matrices,
for rectangular and \cite{CSS08}-type generalized Capelli.
First one is an immediate corollary of holomorphic factorization, 
but surprisingly for rectangular case
we need some additional work. 
We also present \cite{CSS08}-type identity which gives certain conditions
which are sufficient for the Capelli identity.
The first one is a corollary of it, but we think exposition should be organized from simple and explicit to more complicated and abstract.

\subsubsection{Organization of the paper}

Section \pref{SectPrelim} contains some preliminary notations.
We hope reader may skip it and look it only if something is unclear.
In the main text we try to  remind basic notations when we use them
or to make them "intuitively" clear.

Section \pref{SectHolFact} contains the main results.
We prove the holomorphic factorization for the determinants
of certain non-commutative matrices. First we give
some general result, after deduce from it some more concrete propositions
about Capelli matrix.

Section \pref{SectCapel} contains the Capelli identities.
Which we deduce from the results of the previous section.

In section \pref{SectCayl} we mention some analogues
of the Cayley identities which are rather straight-forward.

In the last section \pref{SectLast} we mentions
problems with "dequaternionification" and discuss further open
problems.

\subsubsection{On proofs. Key issue - cancellation is global, but calculation is local \label{SectOnProofs} }

The proofs are quite elementary and easy, however it seems to us that
direct proof of the Capelli identities above would  hardly  be possible,
by the reasons explained below. 
So we hope that (at least part of) simplicity is due to the right strategy.

The main difficulty is proving holomorphic factorization
theorem  \pref{MainTheor} in section \ref{MainThSect}.
All other results follows from it more or less easily.

The key issue about the proof of this result is the following.
The calculational (and so "hard") part of the proof is actually
"local" i.e. we work with $ \Psi \bar \Psi $, but not with the product $\prod_k \Psi_k \bar \Psi_k$.
So calculational part is essentially about 2x2 matrices, which
means it is easily manageable. 
This is proposition  \pref{HolFactPsi}.  
The correction 
$k ->  [ k +1/4 , ~ i/4 ; ~~ i/4 , k - 1/4 ]$ appears here.

So we consider as a key idea  that
it was observed that commutative fact
 $\Psi_{Real} \Psi_{Imag} = \Psi \bar \Psi/(-2i) $
can be relaxed to 
$\Phi_{Real} \Phi_{Imag} = \Psi \bar \Psi/(-2i)  + holomorphic(\Psi_i)$,
where $holomorphic(\Psi_i)$ does not contain $\bar \Psi_i$.
And that this relaxed condition is enough to get desired
results about non-commutative matrices.

The second part is to get from this "local" result
a required theorem. It is quite easy.
See proposition \pref{HolFactGeneral}.
Essentially it is the fact that
$\prod_{k} (\Psi_k \bar \Psi_k/(-2i)  + holomorphic_k(\Psi_i))
= \prod_{k} (\Psi_k \bar \Psi_k/(-2i) )$ 
- too many holomorphic terms multiplies to zero.
I.e. only product $\prod_{k=1...n} \Psi_k \bar \Psi_k$ is non zero,
 but if we have one of terms is not $\Psi \bar \Psi$ but $\Psi_a \Psi_b$
(without "bar") it will be zero.
But the point  is that {\bf unwanted terms cancel only when
we multiply ALL terms $k=1...n$}, for any finite
part $k=1...m<n$ unwanted terms are present and equality above is obviously not true.
So the cancellation is in some sense "global", but "easy".

In this respect the present proof seems to be quite different from
the proofs given in \cite{CFR09} and (probably)
\cite{FZ93}, \cite{CSS08}, where 
proofs go inductively in the size of the matrix and
we get cancellations for each smaller size and use this
on the next step. 
Such strategy is not expected to work for our
Capelli identities, since the cancellation is "global".
The representation theory approach to Capelli identities
\cite{HU91}, \cite{KS91} has a disadvantage that
it cannot work for the general identities of \cite{CSS08}-type,  e.g.
theorem \pref{CCSCap} or our main theorem \pref{MainTheor}.

\subsubsection{Further comments}

Some comments here might be subjective, but to our taste
they might help understanding and  giving some focus for attention.

{\bf What is surprising ?}
The fact that correction changes in quite a non-trivial manner:
\bea
k ->  \left( \begin{array}{ccccccccc}
k +1/4 & i/4  \\
i/4 & k - 1/4
\end{array} \right), ~~~
\eea
seems quite unexpected for us.
When we started we naively  hoped that correction
would be much simpler $k-> diag(k,k)$ or may be $diag(k-1, k+1)$ or may be $1/2$
will appear, but the fact the correction is not real - it contains "i"
and is quite far from integers - contains  $1/4$ seems quite surprising.

{\bf How do we come here ?}
The primary motivation for this paper is \cite{AH11}, where
Capelli-type identity for 4x4 matrices has been obtained.
Since the approach to Capelli identities developed in \cite{CFR09}
seems to us quite simple and powerful, it was interesting
to try to use it on the new case. 
Also after some time of looking on the result of \cite{AH11}
we came to understanding that it can be generalized
to $n\times n$ matrices.
However after some attempts
it became clear that  approach of \cite{CFR09} does not seem to work
directly for the new identities.
Somehow  it become clear that it is much more easy to prove
factorization result and Capelli will follow from it.

{\bf What technique is  used ?} Our approach is rather elementary based
on the algebraic manipulations with anti-commuting (Grassman) variables $\psi_i$.
Generalizing the fact that $\det(M)=\psi^M_1...\psi_n^M$, where
$\psi_i^M= \sum_k \psi_k M_{ki}$. Its usage in non-commutative setup
has been pioneered by Yu. Manin \cite{Manin87}. 
It also has been  used by Japanese researchers in their works
on Capelli identities e.g. \cite{Um98AMS}.
We also applied it \cite{CFR09} to Capelli identities and Manin matrices.

{\bf What seemed "difficult" to us ?}
The first quite simple guess is looking on \cite{AH11} to come to Capelli identities above with unknown form of correction. Second guess
is to understand that it is better to search for holomorphic factorization
result  and reduce Capelli to it. It seems quite important step because 
direct attack on Capelli above is hardly possible - especially on 
anti-symmetric version - which up to now resists all attempts to be 
proved directly by methods like \cite{FZ93},\cite{CSS08},\cite{CFR09},
and also see arguments above that cancellation is somewhat global,
and also for us correction was unknown. 
The third point is somehow find way of finding right correction term
and proving the holomorphic factorization.
In commutative case it can be done due to 
 $\Psi_{Real} \Psi_{Imag} = \Psi \bar \Psi/(-2i) $
So in non-commutative case  taking into account existence of "correction term"
it is natural to search for 
 $(\Psi_{Real}+correction1) (\Psi_{Imag} + corrrection2) = \Psi \bar \Psi/(-2i)$.
 {\bf But this does not exist}. This was main trouble for us.
 At some point working on this problem we lost a hope to get any positive results and just
wanted to get any (since some time already invested) - may be negative result just in the simplest 4x4 case. 
Luckily quite a long direct computations in 4x4 case lead to positive result - we found the identity.
Analyzing it and reassured in happy end,
we see that the only way is to relax condition  $\Psi_{Real} \Psi_{Imag} = \Psi \bar \Psi/(-2i) $
 to $\Phi_{Real} \Phi_{Imag} = \Psi \bar \Psi/(-2i)  + holomorphic(\Psi_i)$.
Things went more or less smoothly after that.

The time spent - some evenings of about 2 months until we get main result and 3 for further improvements, polishing and writing.

\section{Preliminaries \label{SectPrelim} }
The reader is recommended to skip this technical section during the first reading.
For simplicity one may assume that our algebras are over $\RR$
and notation $A^{\CC}$ denotes their complexification: $A\otimes_{\RR} \CC$
(more general formalism is described below).
For a matrix $M$ with elements in $A^{\CC}$
 we  denote by $M^\RR$ - its real form ("decomplexification") 
 i.e. each element $M_{ij}$ is substituted by 2x2 matrix 
$ \left( \begin{array}{ccccccc}
 real(M_{ij}) &  imag(M_{ij}) \\
-imag(M_{ij}) & real(M_{ij})  \\
\end{array} \right) 
 $.  We denote by $\psi^M_k$ the following:  $\psi^M_k=\sum_i \psi_i M_{ik}$ for some variables $\psi_i$.
Column-determinant is defined  by the standard formula $
{\det}^{\mathrm{col}} M=\sum_{\sigma\in S_n} (-1)^\sigma
\prod^{\curvearrowright}_{i=1,...,n} M_{\sigma(i)i}$,  
where $\prod_{i=1,...,n} M_{\sigma(i),i}$ means that the elements
from the first column go first, then from the second column and so on and so
forth.
If necessary see more details below. 


\subsection{Formal complex numbers $\kk[i]$. Notation $A^{\CC}$ } 

For the reader briefly looking this paper we suggest to consider that
our algebras are over $\RR$ and notation $A^{\CC}$ denotes its complexification
$A\otimes_{\RR} \CC$. 
But for those who are interested deeper we provide more general framework.
We can start with an arbitrary ground field $k$ (e.g. finite or $\CC$ ) and 
define its complexification $k[i]$ which  may not be a field but just a commutative ring. The reasons to describe such framework are twofold.
First reason is that clearly considerations here are algebraic
and should work over any (e.g. finite field). Second reason - main examples are related to representation theory so one may want to consider 
representation of  $A^{\CC}$ in some $\CC$-vector space, then we will
face with certain ambiguities - like complex conjugation will act only on 
$A^{\CC}$ or on representation space also ? The formalism below avoids such ambiguities.

Let us introduce a formal variable  $i$ such that  $i^2 = -1$, consider an algebra $\kk[i]$ for some field $\kk$.
We can call $\kk[i]$ "formal complex" numbers.
\footnote{Typically $k= \RR$ and so $k[i]=\CC$, but it is not necessary, in particular $\kk=\CC$ is allowed. More generally $\kk$ can be any commutative ring}

{\bf In this paper we will work with algebras over $\kki$.}

{\Notn ~} For some (may be non-commutative) algebra $A$ over $\kk$, let us  denote its "complexification"  $ A^{\CC} = A+iA = A\otimes_{\kk} \kk[i]$.

{\Notn ~} For any element $ c \in A^{\CC}$ denote by $Re(c)$ its natural projection projection on $A$ (first set in $A+iA$)  and denote by $Im(c) $ projection on the second set.
From now on writing $c = a + ib$ means that $ a= Re(c)$ and $b=Im(c)$. 

{\Notn ~} The map "formal complex conjugation" $c \mapsto \bar c$ is defined as usually $a+ib \mapsto a-ib$.

{\Lem ~} The "formal complex conjugation" $c \mapsto \bar c$ is isomorphism of $A^{\CC}$.

This is obvious.

{\Rem ~} Let us emphasize that even  in the case $\kk = \CC$ all our operations (real, imaginary parts, complex conjugation) does not act on $\kk= \CC$. They take into account only the formal "i". 

%

\subsection{Notation $M^\RR$ - real form of a "complex" matrix $M$ \label{sectRealForms} } 

{\Notn ~}
Consider an algebra $A^{\CC}$ over $\kki$ and some $n$ by $m$ matrix $Z$ with elements in $Z_{ij} \in A^{\CC}$.
Let us call its "decomplexification" (real form)  and denote $Z^{\RR}$ the  $2n$ by $2m$ matrix which is defined
substituting each element $Z_{ij}$ by the two by two matrix 
\bea  \left( \begin{array}{ccccc}
 Z_{R,ij} & Z_{I,ij} \\
 -Z_{I,ij} & Z_{R,ij}
 \end{array} \right) , 
 \eea
where  $ Z_{R,ij} , Z_{I,ij}$ are respectively real and imaginary parts of $Z$.

{\Ex ~}
\bea
Z =  \left( \begin{array}{ccccccc}
 Z_{11} &  Z_{12}  \\
Z_{21} & Z_{22}  \\
\end{array} \right) , ~~~
Z^{\RR} =  \left( \begin{array}{ccccccc}
 Z_{R,11} & Z_{I,11} & Z_{R,12} & Z_{I,12}  \\
 -Z_{I,11} & Z_{R,11} & -Z_{I,12} & Z_{R,12} \\
 Z_{R,21} & Z_{I,21} & Z_{R,22} & Z_{I,22}  \\
 -Z_{I,21} & Z_{R,21} & -Z_{I,22} & Z_{R,22} \\
 \end{array} \right) 
.
\eea  

It is clear that:
{\Prop The correspondence $M\to M^{\RR}$ defines the homomorphism of algebras from $Mat_{n\times n}(A^{\CC}) \to Mat_{2n\times 2n} (A)$.
 }

{\Prop 
Consider elements $\Psi, Z \in A^{\CC}$, write  $\Psi = \Psi_R + i \Psi_I $ and $Z  =  Z_R + i Z_I$,
for real and imaginary parts decomposition, then
\bea
(\Psi_R , \Psi_I)
 \left( \begin{array}{ccccc}
 Z_R & Z_I \\
 -Z_I & Z_R
 \end{array} \right) = ( Re( \Psi Z ), Im (\Psi Z ) ).
\eea  
}
Proof is obvious.
Let us emphasize that nothing more is required in this proposition except correct definitions of $Re, Im$ and $i^2= - 1$. There is no need of commutativity and basic field $\kk$ can be $\CC$. 

Respectively the fact above is true for matrices:
{\Cor \label{CorZRR}
For any elements $\Phi_i, Z_{ij} \in A^{\CC}$ 
\bea
(\Psi_{R,1} , \Psi_{I,1}, \Psi_{R,2} , \Psi_{I, 2}, ..., \Psi_{R,n} , \Psi_{I, n}   )
 \left( \begin{array}{ccccccc}
 Z_{R,11} & Z_{I,11} & Z_{R,12} & Z_{I,12} & ... & Z_{R,1m} & Z_{I,1m} \\
 -Z_{I,11} & Z_{R,11} & -Z_{I,12} & Z_{R,12} & ... & -Z_{I,1m} & Z_{R,1m} \\
 ... & ... &... & ... & ... &... & ... \\
 Z_{R,n1} & Z_{I,n1} & Z_{R,n2} & Z_{I,n2} & ... & Z_{R,nm} & Z_{I,nm} \\
 -Z_{I,n1} & Z_{R,n1} & -Z_{I,n2} & Z_{R,n2} & ... & -Z_{I,nm} & Z_{R,nm} \\ 
 \end{array} \right) 
 = \\ \nn \\  =
 ( Re( \Psi^Z_1 ), Im (\Psi^Z_1 ), Re( \Psi^Z_2 ), Im (\Psi^Z_2 ), ... , Re( \Psi^Z_m ), Im (\Psi^Z_m )  ).
\eea  
where $ \Psi^Z_i$ are (as usually) defined as:
\bea 
(\Psi^Z_1 , \Psi^Z_2 ,  ... , \Psi^Z_m   )
= \\ \nn \\  =
(\Psi_1 , \Psi_2 ,  ... , \Psi_n   )
 \left( \begin{array}{ccccc}
 Z_{11}  & Z_{12}  & ... & Z_{1m} \\
 Z_{21}  & Z_{22}  & ... & Z_{2m} \\
 ...  & ...   & ... & ... \\
 Z_{n1}  & Z_{n2}  & ... & Z_{nm} \\
 \end{array} \right) .
\eea
}

\subsection{Key fact $\Psi_{Real} \Psi_{Imag} = \Psi \bar \Psi/(-2i) $ in commutative case } 

{\Lem \label{PropPsiRIPBP}
Consider $\Psi = \Psi_R + i \Psi_I $,  and respectively $\bar \Psi = \Psi_R - i \Psi_I $,  such that $(\Psi_R)^2=0, (\Psi_I)^2=0, \Psi_R \Psi_I + \Psi_R \Psi_I=0$,
then
\bea
\Psi_R \Psi_I = \frac{1}{ -2i} \Psi \bar \Psi.
\eea
}

Proof is straightforward.

\subsection{Notation $\Psi^M_k$,  Grassman variables,  and column determinant }

{\Notn ~} Consider some $n\times m$ matrix $M$ and some variables $\psi_1, ..., \psi_n$ we will denote by $\psi^M_k$, $k=1...m$ the elements defined as 
$\psi^M_k= \sum_{i=1...n} \psi_i M_{ik}$ or in matrix notations:

\bea ( \psi_1^M, ... ,  \psi_m^M) = (\psi_1, ... ,
\psi_n) \left(\begin{array}{ccc}
 M_{11} & ... &  M_{1m} \\
... \\
 M_{n1} & ... &  M_{nm}
\end{array}\right).
\eea

{\Def \label{coldetdef}
Let $M$ be a matrix with elements in some not necessarily commutative ring.
The column-determinant of $M$  is defined in the standard way with
prescription that one writes at first the elements
from the first column, then from the second column and so on and so
forth: 
\bea
{\det}^{\mathrm{col}} M=\sum_{\sigma\in S_n} (-1)^\sigma
\prod^{\curvearrowright}_{i=1,...,n} M_{\sigma(i)i}, \eea 
where
$S_n$ is the group of
permutations of $n$ letters, and the
symbol $\curvearrowright$ emphasize that in the product
$\prod_{i=1,...,n} M_{\sigma(i),i}$ one writes at first the elements
from the first column, then from the second column and so on and so
forth.}

Consider any algebra $A$ and anticommuting variables $\psi_i \in A$
 (i.e. $\psi_i^2=0, \psi_i \psi_j= -\psi_j
\psi_i$). Consider some elements $M_{pq} \in A$ such that  $\psi_i$ commute with $M_{pq}$: $\forall
i,p,q:~~[\psi_i, M_{pq}]=0$. 


{\Lem \label{detTopFormDef} For an arbitrary matrix $M$ (not
necessarily with commuting entries) it holds: 
\label{DetTopForm}
\bea
 {\det}^{\mathrm{col}} (M) ~ \psi_1  ...  \psi_n
= \mpsi_1... \mpsi_n, \eea  } 
The proof is straightforward.

Let us emphasize the $\mpsi_i$ will not in general anticommute. The requierment
of their anticommutativity is one of the definitions of \MM \cite{CFR09},
where one can find further discussions of non-commutative determinants, etc.

\section{Holomorphic factorization of the determinants \label{SectHolFact}  }

This section contains the results on holomorphic factorization of the determinants
of the form $det(M^\RR+cor_R)= det(M+cor)det(\bar M+ \widebar{cor} )$. 
Here are 3 subsections 1) we remind and analyze the proof of this fact in the commutative case,
obtaining the first (easy) non-commutative generalization for matrices where elements in each column
commute (in particular for Manin matrices) 2) we prove our main result - theorem \pref{MainTheor}
which describes certain quite general commutation relations which ensures factorization theorem
3) we prove factorization theorem for Capelli type matrix - it rather easily follows from our main theorem.

\subsection{Holomorphic factorization in commutative and Manin cases}
 
Let us recall the standard proposition $det(M^\RR)= det(M)det(M)$ for matrices over commutative rings and analyze its proof (looking how much commutativity is used). This proof and analysis will be guiding principle for what follows,
we strongly suggest the reader who is interested in the proof of the main theorem not to omit this section.
As a bonus of our analysis we prove first non-commutative (easy) version of the factorization theorem.
This easy version does not contain "corrections" (i.e. $det(M)$, not $det(M+{cor} )$).

{\Prop  \label{HolFactComProp}
Consider an $n$ by $n$  matrix  $C$ over commutative alegbra $A^{\CC}$.
Let $2n$ by $2n$ matrix $C^{\RR}$ with elements in $A$  be its "decomplexification" (real form),  \footnote{ if necessary - see subsection \ref{sectRealForms} for notations} 
then
\bea
 {\det}^{\mathrm{ }}_{2n \times 2n}  (C^{\RR} ) =
{\det}^{\mathrm{}}_{n \times n}  (C )  {{\det}^{\mathrm{}}_{n \times n}  (\bar C ) }.  
\eea
}

\PRF
Consider Grassman algebra $\Lambda[\psi_{R,1}, \psi_{I,1},..., \psi_{R,n}, \psi_{I,n}]$ with $2n$  generators 
, define $n$ variables $\psi_i =   \psi_{R,1}+ i \psi_{I,1}  $ 
which are elements of the complexified Grassman algebra $\Lambda^{\CC}[\psi_{R,i}, \psi_{I,i} ] $.  

The proof goes as follows: 
\bea
{\det}^{\mathrm{}}_{2n \times 2n}  (C^{\RR} ) \prod_{k=1...n} \psi_{R,k} \psi_{I,k} = \\
\mbox{ use lemma  \ref{DetTopForm} page \pageref{PropPsiRIPBP}}  \\ 
%
=\prod_{k=1...2n} \psi^{C^{\RR}}_{k}  
= \\
\mbox{ use lemma \ref{PropPsiRIPBP} page  \pageref{PropPsiRIPBP}} \\
= \prod_{k=1...n} Re(\psi^{C}_k) Im({ \psi }^{C}_k) 
= \label{repsiimpsiInProof1} \\
\mbox{ use corollary \ref{CorZRR} that } \psi^{Real} \psi^{Imag} = \psi \bar \psi /(-2i)  \\
= \frac{1}{(-2i)^n }\prod_{k=1...n} \psi^{C}_k \bar{ \psi }^{C}_k 
= \label{holcomantiInProof1} \\
\mbox{ use that }  \psi \bar \phi = - \bar \phi \psi  \\
= \frac{1}{(-2i)^n }  (-1)^{n(n-1)/2} 
\prod_{k=1...n} \psi^{C}_k \prod_{k=1...n}  \bar{ \psi }^{C}_k  =  \\
\mbox{ use lemma \ref{DetTopForm} page \pageref{DetTopForm} }  \\ 
= \frac{1}{(-2i)^n } (-1)^{n(n-1)/2}
{\det}^{\mathrm{}}_{n \times n}  (C ) \prod_{k=1...n} \psi_k 
~~ { {\det}^{\mathrm{}}_{n \times n}  
(\bar C )   }  \prod_{k=1...n}  \bar{ \psi }_k  = \\
\mbox{ use   } \psi^{Real} \psi^{Imag} = \psi \bar \psi /(-2i) \mbox{ use   }  \psi \bar \phi = - \bar \phi \psi \\ 
= 
{\det}^{\mathrm{}}_{n \times n}  (C )  { {\det}^{\mathrm{} } _{n \times n}  
(\bar C )   }  \prod_{k=1...n} \psi_{R,k} \psi_{I,k}.
\eea
Comparing the first and the last expressions we get the desired result.

Proposition is proved. \BX

Let us analyze what properties has been used in this proof.

We have used lemma  \ref{DetTopForm}, which
states:
\bea
%
{\det}^{\mathrm{col}}_{n \times n}  (M ) \prod_{k=1...n} \psi_{k} 
=\prod_{k=1...n} \psi^{M}_{k}.  \label{DetViaGrassInAnalyse}
\eea
This lemma does not require any commutativity or other relation for elements of  matrix $M$
if we use column-determinant. \footnote{See definition \pref{coldetdef}. Briefly speaking 
 it is defined by the same rule as the standard one, with the convention that
elements from first column are at the first position in the products, from second column are at the second position, etc.} Clearly for matrices with commuting entries column-determinant is standard determinant.

We used:
\bea
\psi \bar \phi = - \bar \phi \psi .
\eea
This is unavoidable constraint in our approach and we will always require 
that we work with some elements $M_i$, such that $[M_i, \bar M_j] = 0$.
This is enough for our main application to the Capelli identities.

And finally the last requirement:
\bea
\psi^{Real} \psi^{Imag} = \frac{1}{-2i} \psi \bar \psi.
\eea
Actually the main soul of the paper is around this formula. It is impossible to get the direct analog of it in interesting cases, but some modification is necessary and  the main point is that terms deviating from this commutative formula can be
finally canceled.

%
%
%
%

The following simple statement is our first bonus for careful looking on 
relaxation of commutativity.

{\Th \label{theor1} Consider an $n$ by $n$  matrix $M$ over $A^{\CC}$ such	 that
for any two its matrix elements 
$\alpha, \beta$, it is true that $[\alpha, \bar \beta]=0$, i.e.
all elements commute with complex conjugated elements; 
and elements in each column commute among themselves: $\forall ~i,j, c ~ [M_{ic}, M_{jc} ] = 0$ (for example this is true if $M$ is a \MM~ (see subsection \pref{SSSectManinMatrReminder} or \cite{CFR09}) ), then holomorphic factorization for the determinant 
holds true :
\bea
{\det}^{\mathrm{col}}_{2n\times 2n} M^{\RR} = {\det}^{\mathrm{col}}_{n\times n} M  ~{\det}^{\mathrm{col}}_{n\times n} \bar M .
\eea
}

Here ${\det}^{\mathrm{col}}$ is column-determinant.
{See definition \pref{coldetdef}. }

\PRF
The proof above actually works in this more general case.

Indeed, the condition that elements from the same column commute ensures that 
\bea
Re(\psi^{M}_k) Im({ \psi }^{M}_k) 
= \frac{1}{(-2i) } \psi^{M}_k \bar{ \psi }^{M}_k ,
\eea
and hence equality \ref{repsiimpsiInProof1} holds true.

Condition that $[\alpha, \bar \beta]=0$ for all matrix elements  $\alpha, \beta$ of $M$,
 ensures that equality \ref{holcomantiInProof1} is also true.
 
As we discussed above  equality \ref{DetViaGrassInAnalyse}
holds without any conditions.

So we see that all steps of the proof of the proposition \ref{HolFactComProp}
also work in the more general setup of the theorem above.

Theorem is proved.
\BX

\subsection{Main theorem \label{MainThSect}}
This subsection contains the main theorem of this paper - we prove that under certain
quite general conditions on matrix with non-commutative elements we can prove
factorization result $det(M^\RR+cor_R)= det(M+cor)det(\bar M+ \widebar{cor} )$. 
The proof is divided into 3 steps. We have discussed the key ideas of the proof in
section \pref{SectOnProofs}. Also for those who are interested in proof we strongly 
recommend first to look in the previous subsection where toy model is analysed.

Consider some Grassman algebra $\Lambda[\psi_i]$. Let us recall that for any matrix $M$ we denote by  $\psi^M_l$ expressions defined as $\psi^M_l = \sum_k \psi_k M_{kl}$.

{\Th \label{MainTheor} Consider matrices $C,Q$ with elements in $A^{\CC}$. Assume that
for any two matrix elements $\alpha, \beta$ of matrices $C,Q$, it is true that $[\alpha, \bar \beta]=0$, i.e.
all elements commute with complex conjugated elements; assume also:
\bea
\forall k, ~~(\psi^C_k)^2 = \psi^C_k \psi^Q_k, ~~  \psi^C_k  \psi^Q_k  = - \psi^Q_k \psi^C_k ,
~~
(\psi^Q_k)^2  = 0, 
\eea
(Or the same requirements can be written more explicitly in terms of matrix elements):
\bea
\forall i,j,k: ~~[C_{ik}, C_{jk}] = C_{ik} Q_{jk} -  C_{jk}Q_{ik},  
~~~ [C_{ik}, Q_{jk}] = [C_{jk}, Q_{ik}],
~~~ [Q_{ik}, Q_{jk}] = 0, 
\eea
pay attention that  index $k$ is the same for both $\psi^C$, $\psi^Q$, i.e. all
relations concern elements from the same ($k$-th) column, there is NO relations between elements from different columns,

{\bf then} the following analog of holomorphic factorization holds:

for any central elements $d_i$:
\bea
{\det}^{\mathrm{col}}_{2n\times 2n} (C^{\RR} + Q^{\RR} CorrTriDiag ) = \\ =
{\det}^{\mathrm{col}}_{n\times n} (C + Q diag(d_n, d_{n-1}, ...,d_{2}, d_{1} ) ) 
{\det}^{\mathrm{col}}_{n\times n} (\bar C +\bar Q diag( d_n,  d_{n-1}, ...,  d_{2}, d_{1} ) ), 
\eea
where there are two options for the matrix $CorrTriDiag$ , 
the first one is:
\bea
CorrTriDiag  = \left( \begin{array}{ccccccccc}
d_n +1/4 & i/4 & 0 & 0 & 0 & 0 & ... & 0 & 0 \\
i/4 & d_n - 1/4 & 0 & 0 & 0 & 0 & ... & 0 & 0 \\
 0 & 0 & d_{n-1} +1/4 & i/4 & 0 & 0 & ... & 0 & 0 \\
 0 & 0 & i/4 & d_{n-1} - 1/4 &0 & 0 & ... & 0 & 0 \\
 ... & ... & ... & ... & ...  & ... & ... & ... & ... \\
 0 & 0 & 0 & 0 & 0 & 0 & ... & d_{1}+1/4 & i/4 \\
 0 & 0 & 0 & 0 &0 & 0 & ... & i/4 & d_1 -1/4 
\end{array} \right). ~~~
\eea
The second option for $CorrTriDiag$ is the same formula with substitution $i$ to $-i$,
i.e. complex conjugate to the first one. 
}

Here ${\det}^{\mathrm{col}}$ is column-determinant.
{See definition \pref{coldetdef}. }

{\Rem ~} The elements $Q_{ik}$ from the same column $k$, $i=1...n$ form a commutative ring. If we strengthen  the second relation to the form
$0= [C_{ik}, Q_{jk}] = [C_{jk} Q_{ik}]$, then 
elements $C_{ik}$ form a universal enveloping of Lie algebra over this
ring, moreover each pair $C_{ik}, C_{jk}$ is Lie subalgebra.
But without such strengthening we think that the algebra generated by 
 $C_{ik},Q_{jk}$ would have much bigger size, we think it is of exponential
growth. It would be interesting to calculate the Poincare series
of this algebra.\footnote{ \label{thankMO}
It is easy to see that the relations above imply that $(\psi_k^C)^3=0$, because
$(\psi_k^C)^3=(\psi_k^C) (\psi_k^C)^2= (\psi_k^C) (\psi_k^C\psi_k^Q)
=(\psi_k^C)^2 (\psi_k^Q)=(\psi_k^C \psi_k^Q) (\psi_k^Q)= 0 $,
the Poincare series of the simplest algebra with 3 generators $a_i$ and
relation $(\sum_i a_i \psi_i)^3 =0$ 
has been discussed at Mathoverflow \href{http://mathoverflow.net/questions/81415}{mathoverflow.net/questions/81415}
\href{http://mathoverflow.net/questions/82722}{mathoverflow.net/questions/82722}  
.
 David Speyer proposed the answer $\frac{1}{1-3x+x^3}$ for the algebra
 with 3 generators $a_i$ and relation $(\psi)^3=0$ with $\psi=\sum_{i=1,2,3} a_i\psi_i$, Vladimir Dotsenko proposed use of Anick's resolution to explain the answer.}

 
{\Rem ~} The relations above are generalizations of the relations for matrix
elements of Capelli matrix $E_{ij}$, initially we aimed to prove the theorem
for it, but when the proof has been obtained we see that it works
in generality above. In Capelli case $Q_{ij} = \delta_{ij}$ and there are
many relations for elements from  different columns which are unnecessary
for this theorem. 

{\Rem ~} Any submatrix of matrix $C$ is again matrix satisfying the conditions
of the theorem with the matrix $\tilde Q$ which is corresponding submatrix of $Q$.

\subsubsection{Example: n= 1 \label{SSectN1} }

Let us consider the example of the main theorem for $n=1$, which corresponds to $2 \times 2$
decomplexified matrices. The case $n=2$ leads to $4\times 4$ matrices and 
it is seems such an example will be too long.  It seems the shortest way 
to see the result in $n=2$ case is to repeat the proof. 

The case $n=1$ is actually quite specific and does not demonstrate all features of the general case:
we need to impose less conditions than in general case and moreover the choice of correcting matrix 
$CorrTriDiag$ has big degrees of freedom, in particular, it can be taken to be zero, while
already for $n=2$ we will see that the correcting matrix is fixed quite uniquely.
Still it might be useful to present this example.

The following proposition illustrates the main  theorem  for $n=1$.

{\Prop \label{n1pr1} Consider $C,Q$, such that $0= [C, \bar  C] =[C, \bar  Q ] = [Q, \bar  C] =[Q, \bar  Q ]$,
i.e. all complex-conjugated elements commute with $C,Q$.  Consider also some central element $d$,
then:
\bea
{\det}^{\mathrm{col}}_{2\times 2}
\left( 
\left( \begin{array}{ccccccccc}
Re(C) & Im(C) \\
-Im(C) & Re(C) \\
\end{array} \right)
+
\left( \begin{array}{ccccccccc}
Re(Q) & Im(Q) \\
-Im(Q) & Re(Q) \\
\end{array} \right)
\left( \begin{array}{ccccccccc}
d+1/4 & i/4 \\
i/4 & d-1/4 \\
\end{array} \right)
\right) 
= \\
= (C+ d Q)( \bar C + d \bar Q).
\eea
}

{\Rem ~} It might be surprising that if $d=0$, then the left hand side contains $Q$, while the right hand side does not - we will  see this is indeed correct. Comparing to the main theorem we see that for $n=1$
we do not need the condition $[C,Q]=0$.
 
Before proving this proposition let us give another which illustrates the mechanism why it works:

{\Prop \label{n1pr2} Consider $A = a + ib, Q= p+iq $, such that $[Q, \bar A]=0$, then:
\bea \label{prop222}
{\det}^{\mathrm{col}}_{2\times 2}
\left( 
\left( \begin{array}{ccccccccc}
a & b \\
- b  & a \\
\end{array} \right)
+
\left( \begin{array}{ccccccccc}
p & q \\
-q & p \\
\end{array} \right)
\left( \begin{array}{ccccccccc}
1/4 & i/4 \\
i/4 & -1/4 \\
\end{array} \right)
\right) 
=  a^2+b^2.
\eea

If we additionally assume that $[A, \bar A]=0$, then it also equals to:
\bea
= A \bar A.
\eea
}
{\Rem ~ \label{EigRem} } Pay attention that there is no $Q$ at the right hand side. 
To comment on this somewhat surprising fact, let us mention that it is s probably related to the fact,
that the matrix 
$\left( \begin{array}{ccccccccc}
1/4 & i/4 \\
i/4 & -1/4 \\
\end{array} \right)
$ contains eigenvectors of any decomplexified 2x2 matrix, i.e. matrix of the form:
$\left( \begin{array}{ccccccccc}
p & q \\
-q & p \\
\end{array} \right)$.

\PRF  Let us prove proposition \ref{n1pr2}.
\bea
\left( \begin{array}{ccccccccc}
p & q \\
-q & p \\
\end{array} \right)
\left( \begin{array}{ccccccccc}
1/4 & i/4 \\
i/4 & -1/4 \\
\end{array} \right)
=1/4\left( \begin{array}{ccccccccc}
p + iq & ip-q\\
-q+i p  & -p - i q \\
\end{array} \right)
=(p+iq)/4\left( \begin{array}{ccccccccc}
1  & i  \\
i  & -1 \\ 
\end{array} \right)
=Q/4\left( \begin{array}{ccccccccc}
1  & i  \\
i  & -1 \\
\end{array} \right)
\eea
Hence the left hand side in the proposition can be transformed:
\bea
{\det}^{\mathrm{col}}_{2\times 2}
\left( 
\left( \begin{array}{ccccccccc}
a & b \\
- b  & a \\
\end{array} \right)
+
Q/4 \left( \begin{array}{ccccccccc}
1 & i \\
i & -1 \\
\end{array} \right)
\right) 
= 
{\det}^{\mathrm{col}}_{2\times 2}
\left( \begin{array}{ccccccccc}
a+Q/4 & b+iQ/4 \\
- b+iQ/4  & a-Q/4 \\
\end{array} \right)
= \\
(a+Q/4)(a-Q/4)- (-b +iQ/4) (b+iQ/4)
= \\ =
a^2 - aQ/4+Q/4a- (Q/4)^2 +b^2 +biQ/4-iQ/4b + (Q/4)^2
= \\ =
a^2 + b^2 +[Q/4,a]-[Q/4, ib]
=
a^2 + b^2 +[Q/4,a-ib]
=
a^2 + b^2 +[Q/4, \bar A].
\eea
Now we can use that $[Q, \bar A]=0$ and so we finally get that the left hand side of  \pref{prop222} equals to:
\bea
= a^2 + b^2.
\eea

So we derived the first claim of the proposition. The second claim follows trivially indeed:

$[A,\bar A]=0$ implies that $[a, b]=0$ and hence $A \bar A =a^2 + b^2$.

\BX 
 
\PRF  Let us prove proposition \ref{n1pr1}.

\bea
{\det}^{\mathrm{col}}_{2\times 2}
\left( 
\left( \begin{array}{ccccccccc}
Re(C) & Im(C) \\
-Im(C) & Re(C) \\
\end{array} \right)
+
\left( \begin{array}{ccccccccc}
Re(Q) & Im(Q) \\
-Im(Q) & Re(Q) \\
\end{array} \right)
\left( \begin{array}{ccccccccc}
d+1/4 & i/4 \\
i/4 & d-1/4 \\
\end{array} \right)
\right) 
= \\ =
{\det}^{\mathrm{col}}_{2\times 2}
\left( 
\left( \begin{array}{ccccccccc}
Re(C) +d Re(Q) & Im(C) +d Re(C) \\
-Im(C) - d Im(Q) & Re(C) + d Re(Q) \\
\end{array} \right)
+
\left( \begin{array}{ccccccccc}
Re(Q) & Im(Q) \\
-Im(Q) & Re(Q) \\
\end{array} \right)
\left( \begin{array}{ccccccccc}
1/4 & i/4 \\
i/4 & 1/4 \\
\end{array} \right)
\right).
\eea
So we see  that we can apply the proposition \ref{n1pr1} and get the desired result.
\BX
{\Rem ~}
It may seems that we need to  require $d$ to be real, however more carefully looking we can conclude it is not really necessary. For example the right hand side can be  written as:
\bea
 (C+ d Q)( \bar C + d \bar Q) = (Re(C) + d Re(Q) + i Im(C) + idIm(Q) ) ( Re(C) + d Re(Q) - i Im(C) - idIm(Q) ) = \\
 = (Re(C) + d Re(Q)) ^2 +  ( Im(C) + dIm(Q) ) ^2 + 
 - i [ Re(C) + d Re(Q) ,  Im(C) + dIm(Q) ] =
\eea
We should note that $[ Re(C) + d Re(Q) ,  Im(C) + dIm(Q) ] =
1/4i [ C+dQ + \bar C + d \bar Q, C+dQ - \bar C - d \bar Q] =0$, since bar elements commute with non-bar
and $d$ is central. So we get that the left hand side is equal to: 
\bea
= (Re(C) + d Re(Q)) ^2 +  ( Im(C) + dIm(Q) ) ^2.
\eea

\subsubsection{Proof Step 1. Factorization of $2n\times 2n$ and generalization of
Generalization of $-2i\psi_R  \psi_I  = \psi \bar \psi $}

Consider Grassman algebra $\Lambda[\psi_{R,1}, \psi_{I,1},..., \psi_{R,n}, \psi_{I,n}]$ with $2n$  generators, 
define $n$ variables $\psi_k =   \psi_{R,k}+ i \psi_{I,k}  $ 
which are elements of the complexified Grassman algebra $\Lambda^{\CC}[\psi_{R,1}, \psi_{I,1},..., \psi_{R,n}, \psi_{I,n} ] $.
For $n\times n$ matrix $C$ as usually we denote by
$\psi_k^C=\sum_i \psi_i C_{i,k}$, 
for $2n\times 2n$ matrix $R$ we use 
similar notation:  
$\psi_{k}^R  = \sum_{i=1...n} \psi_{R,i}R_{2i-1,k} + \psi_{I,i}R_{2i,k}$,
where $k=1...n$.

{\Prop \label{HolFactGeneral}
Consider any $2n\times 2n$ matrix $R$ with elements in any algebra $B$, consider $n\times n$ matrices  $C_1,C_2$, with mutually commuting entries (i.e. $\forall k,l,p,q: ~[(C_{1})_{kl},(C_{2})_{pq} ] =0 $).  
Then

in order to have the equality:
\bea
{\det}^{\mathrm{col}}_{2n\times 2n} (R) 
={\det}^{\mathrm{col}}_{n\times n} (C_1) ~
{\det}^{\mathrm{col}}_{n\times n} (C_2 ) , 
\eea

it is sufficient that there exists some elements $b_{k,pq}\in B$, such that the following is true for all $k$:
\bea
\psi^{R}_{2k-1} \psi^{R}_{2k}  =  \frac{1}{-2i} \psi^{C_1}_k   {\bar \psi}^{C_2}_k   + \sum_{pq} b_{k,pq} \psi_p \psi_q . 
\eea
I.e. purely "holomorphic" term $\sum_{pq} b_{k,pq} \psi_p \psi_q $ (i.e. term
not containing $\bar \psi$) is allowed to correct commutative formula \pref{fmlHolFactInMainSect} without breaking factorization property.   

Another sufficient condition is to put purely antiholomorphic term:
\bea
\psi^{R}_{2k-1} \psi^{R}_{2k}  =  \frac{1}{-2i} \psi^{C_1}_k   {\bar \psi}^{C_2}_k   + \sum_{pq} b_{k,pq} \bar \psi_p \bar \psi_q . 
\eea
}

\PRF
Indeed  
when we multiply  over $k$  the terms
on the left hand side and using lemma \pref{DetTopForm}  we get: 
\bea
{\det}^{\mathrm{col}}_{2n \times 2n}  (R ) \prod_{k=1...n} \psi_{R,k} \psi_{I,k}.
\eea

Consider the multiplication of the terms on the right hand side.
Key (but simple) fact is the following:
\bea \label{keyFactSimple} 
\prod_{k=1...n} 
(\frac{1}{-2i} \psi^{C_1}_k   {\bar \psi}^{C_2}_k   + \sum_{pq} b_{k,pq} \psi_p \psi_q)=
\prod_{k=1...n} 
\frac{1}{-2i} \psi^{C_1}_k   {\bar \psi}^{C_2}_k =
\eea
Let us call terms containing only $\psi_k$, but not $\bar \psi_l$ by "holomorphic".
This is true because the product which contains $l$ terms of the type
$\sum_{pq} b_{k,pq}  \psi_p  \psi_q$ and $n-l$ terms of type 
$(\psi^{C_1}_k  { \psi}^{C_2}_k    )$
will have $2n-2l$ holomorphic generators $\psi_t$ from the first one and $l$ from the second type, in total $2n-l=n+(n-l)$, but any product containing more that $n$ holomorphic $\phi$ is zero.  
So we proved this identity.

Let us mention that it is obviously not true if we restrict the product 
to any finite $k=1...m<n$.

Now using that $\psi^{C_1}_r $ and $\bar \psi^{C_2}_{t}$ anticommute
we can continue:
\bea
=(\frac{1}{-2i})^k (-1)^{n(n-1)/2} \prod_{k=1...n} \psi^{C_1}_k  
    \prod_{k=1...n} {\bar \psi}^{C_2}_k =
\eea
and applying lemma \pref{DetTopForm}  twice
we can be continue:
\bea
=(\frac{1}{-2i})^n (-1)^{n(n-1)/2} 
{\det}^{\mathrm{col}}_{n\times n} (C_1) \prod_{k=1...n} \psi_k  ~~
{\det}^{\mathrm{col}}_{n\times n} (C_2) \prod_{k=1...n} \bar \psi_k =   
\eea
using  that $\psi_{R,k}\psi_{I,k}=(1/(-2i))\psi_k \bar \psi_k$ ( $\psi_k$ do not contain any non-commuting variables so it is just standard simple fact) we get:	 
\bea \label{LocFml1600}
{\det}^{\mathrm{col}}_{n\times n} (C_1) ~
{\det}^{\mathrm{col}}_{n\times n} (C_2)  \prod_{k=1...n} \psi_{R,k}\psi_{I,k}.
\eea

Comparing expressions in front of 
$\prod_{k=1...n} \psi_{R,k} \psi_{I,k}$ in formulas \ref{keyFactSimple} and \ref{LocFml1600}
we get the desired result.

Proof of the proposition \ref{HolFactGeneral} finished.
\BX

\subsubsection{Proof Step 2.  Calculation of the correction}

In commutative case:
\bea
\label{fmlHolFactInMainSect}
Re(\psi^{C}_k) Im(\psi^{C}_k)  =  \frac{1}{-2i} (\psi^C_k) ( {\bar \psi^C_k} )	 . 
\eea

In our notations:
\bea
Re(\psi^{C}_k) = \psi^{C^\RR}_{2k-1}, ~~~ Im(\psi^{C}_k) = \psi^{C^\RR}_{2k}
\eea

Our substitute in non-commutative case is the following:

{\Prop \label{HolFactPsi} Assume that $\phi, \psi \in A^{\CC} $ are such that $\bar \phi, \bar \psi$
anticommutate with  them  (but anticommutativity between $\phi,\psi$ is not assumed).
Consider some element $k$ which commute with $\phi, \psi, \bar \phi, \bar \psi$,
then:

in order factorization
\bea
(-2i) (Re(\psi) +  a/2 \phi + b/2 \bar \phi)( Im(\psi)  + c/2i \phi + d/2i \bar \phi)
=  (\psi + k \phi) ( \bar \psi + k \bar \phi ) + hol(\phi,\psi) 
\eea
to exist, (where $hol(\phi,\psi)$ is expression depending on $\phi, \psi$, but not $\bar \phi, \bar \psi$),
 it is sufficient that $a,b,c,d$ satisfy either the conditions:
\bea \label{HolFactPsiFistCond}
a = k,~~ c = k,~~ d= b- 2k, ~~ (\bar \psi + b \bar \phi) ( -\bar \psi + (b -2k) \bar \phi) =0;  
\eea
or conditions:
\bea
b = k,~~ d = -k,~~ c= 2k -a, ~~ (\bar \psi + k \bar \phi)^2 = 0.
\eea
Moreover these conditions are necessary if pair products $ \phi \bar \psi, \psi \bar \phi, \phi \bar \phi,
(\psi +a \phi) (\psi + c \phi), (\bar \psi +b \bar \phi) (-\bar \psi + d \bar \phi) $ are linearly independent
and $a,b,c,d,k$ belong to ring without zero divisors.
}

{
\PRF
Let us transform left hand side:
\bea
(-2i) (Re(\psi) +  a/2 \phi + b/2 \bar \phi)( Im(\psi)  + c/2i \phi + d/2i \bar \phi)
= \\ =
(-2i)/(2*2i) (\psi + \bar \psi  +  a \phi + b \bar \phi)( \psi - \bar \psi  + c \phi + d \bar \phi) 
= \\ =
\frac{1}{-2}( 
 (\psi +   a \phi  )( \psi + c \phi ) 
+ 
 ( \bar \psi  +   b \bar \phi)(  - \bar \psi   + d \bar \phi) 
+
 ( \psi + a \phi )(  - \bar \psi  + d \bar \phi) 
+  (  \bar \psi  + b \bar \phi)( \psi + c \phi )  )
= \\ =
\frac{1}{-2}( 
 (\psi +   a \phi  )( \psi + c \phi ) +  ~~~~ \mbox{  (holomorphic term) }  \label{holTermLHS} \\  
+ 
 ( \bar \psi  +   b \bar \phi)(  - \bar \psi   + d \bar \phi)+ ~~~~ \mbox{ (antiholomorphic term) }  \label{antiHolterm} \\ 
+
 (-2)\psi \bar \psi -(a+c) \phi \bar \psi +(d-b) \psi \bar \phi +(ad - bc) \phi \bar \phi  \label{holAntiholLHS}
 )
\eea
On the other hand, right hand side is:
\bea
(\psi + k \phi) ( \bar \psi + k \bar \phi ) = 
\psi  \bar \psi + k \phi \bar \psi + k \bar \phi \psi + k^2 \phi  \bar \phi  \label{holAntiholRHS}
\eea
In order to have an equality, we should require:
\bea
-(a+c) = -2k, ~~ d-b = -2k, ~~ ad-bc = -2k^2, ~~( \bar \psi  +   b \bar \phi)(  - \bar \psi   + d \bar \phi) =0.
\eea
Since: the first three ensures \ref{holAntiholLHS} equals to \ref{holAntiholRHS},
the last equality ensures the antiholomorphic term \ref{antiHolterm} vanishes.
So the difference between left hand side and right side becomes only holomorphic term \ref{holTermLHS},
which is  $hol(\phi, \psi)$.

The first three can be transformed:
\bea
c= 2k -a, ~~ d = b -2k, ~~ (a-k)(b-k)=0.
\eea
Assuming absence of zero divisors we come to two possible options:
first one with $a=k$:
\bea
a=k, c= k, ~~ d = b -2k, ~~ ( \bar \psi  +   b \bar \phi)(  - \bar \psi   + d \bar \phi) = 0,
\eea
and the second one with $b=k$:
\bea
b=k, d= -k, ~~ c = 2k-b, ~~ ( \bar \psi  +   k \bar \phi)(  - \bar \psi   + (-k) \bar \phi) = 0.
\eea
So the  proposition \ref{HolFactPsi} is proved.
\BX
}
 
Similarly we can prove, factorization modula anti-holomorphic terms:

{\Prop \label{HolFactPsimodAnti} Assume that $\phi, \psi \in A^{\CC} $ are such that $\bar \phi, \bar \psi$
anticommutate with  them  (but anticommutativity between $\phi,\psi$ is not assumed).
Consider some element $k$ which commute with $\phi, \psi, \bar \phi, \bar \psi$,
then:

in order factorization
\bea
(-2i) (Re(\psi) +  a/2 \phi + b/2 \bar \phi)( Im(\psi)  + c/2i \phi + d/2i \bar \phi)
=  (\psi + k \phi) ( \bar \psi + k \bar \phi ) + antihol(\bar \phi, \bar \psi) 
\eea
to exist, (where $antihol(\bar \phi,\bar \psi)$ is expression depending on $\bar \phi, \bar \psi$, but not $\phi,  \psi$),
 it is sufficient that $a,b,c,d$ satisfy either the conditions:
\bea
a = k,~~ c = k,~~ d= b- 2k, ~~ (\psi + k \phi)^2 =0;  
\eea
or conditions:
\bea \label{cond2modAnti}
b = k,~~ d = -k,~~ c= 2k -a, ~~ (\psi + a  \phi)(\psi + (2k-a)  \phi) = 0.
\eea
} 

The propositions above were quite general, not related to the particular
setup of the main theorem we are proving. Now  let us recall the setup and 
show that the conditions in the theorem guarantee that 
we can apply propositions above. And write down explicitly  resulting expressions.

{
\Cor \label{CorOnFact} Assume same as in previous proposition that $\phi, \psi$ anticommute with $\bar \phi, \bar \psi$,
and additionally $\psi^2 = \psi \phi, \phi \psi = - \psi \phi, \phi^2=0$, then 
\bea \label{CorOnFactFmlLHS}
(Re(\psi) + (k+1/4)Re(\phi) + (i/4) Im(\phi) )  (Im(\psi) + (+i/4)Re(\phi) + (k-1/4) Im(\phi) )  
= \\ =
\frac{1}{-2i}(\psi + k \phi )  ( \bar \psi + k \bar \phi  )  + antihol(\bar \psi,\bar \phi)  ,
\eea 
where $antihol(\bar\psi,\bar \phi) = (\phi + k \psi)^2$ - it depends only on $\phi, \psi$.

And also similar equality, obtained by conjugation:
\bea
(Re(\psi) + (k+1/4)Re(\phi) - (i/4) Im(\phi) )  (Im(\psi) + (-i/4)Re(\phi) + (k-1/4) Im(\phi) )  
= \\ =
\frac{1}{-2i}(\psi + k \phi )  ( \bar \psi + k \bar \phi  )  + antihol(\bar \psi, \bar \phi).  
\eea  
}

{\PRF
Consider conditions \ref{cond2modAnti}: $b = k,~~ d = -k,~~ c= 2k -a, ~~ (\psi + a  \phi)(\psi + (2k-a)  \phi) = 0$.
Let us transform:
\bea
(\psi + a  \phi)(\psi + (2k-a)  \phi) = 
\psi^2 + a\phi \psi + (2k-a)\psi \phi + a(2k-a) \phi^2 
= 
\psi \phi  + (2k-2a)\psi \phi  
\eea
So we must require: $1= 2a-2k$ , so $a=k+1/2$, hence $c=k-1/2$.

Hence:
\bea
(-2i) (Re(\psi) +  a/2 \phi + b/2 \bar \phi)( Im(\psi)  + c/2i \phi + d/2i \bar \phi)
= \\ =
(-2i) (Re(\psi) +  (k/2+1/4) \phi + k/2 \bar \phi)( Im(\psi)  + (k-1/2) /2i \phi + (-k)/2i \bar \phi)
= \\ =
(-2i) (Re(\psi) +  (k+1/4) Re(\phi) + (i/4) Im(\phi) )( Im(\psi)  + (-1/4i)  Re(\phi) + (2k-1/2)i/2i Im( \phi) )
= \\ =
(-2i) (Re(\psi) +  (k+1/4) Re(\phi) + (i/4) Im(\phi) )( Im(\psi)  + (i/4)  Re(\phi) + (k-1/4) Im(\phi) ),
\eea
so we get  \ref{CorOnFactFmlLHS}, and so the first claim of the corollary is proved.

To get the second claim we can apply conjugation,
either deduce it in the same way from \ref{HolFactPsiFistCond}.

Corollary \ref{CorOnFact} is proved.
\BX
}

\subsubsection{ Proof Step 3. Finishing proof}

Now the proof of the main theorem \ref{MainTheor} follows easily from proposition \ref{HolFactGeneral} and corollary \ref{CorOnFact}.

\PRF
Consider $\psi = \psi^C_k$, $\phi = \psi^Q_k$ and $k=d_k$. They satisfy the conditions of the corollary \ref{CorOnFact}.

Hence the factorization exists:
\bea
(Re(\psi^C_k) + (d_k+1/4)Re(\psi^Q_k) + (i/4) Im(\psi^Q_k) ) *
\\  * (Im(\psi^C_k) + (+i/4)Re(\psi^Q_k) + (d_k-1/4) Im(\psi^Q_k) )  
= \\ =
\frac{1}{-2i}(\psi^C_k + d_k \psi^Q_k )  ( \bar \psi^C_k + d_k \bar \psi^Q_k  )  + antihol(\bar \psi^C_k,\bar \psi^Q_k)  .
\eea 
 
So we can apply the proposition \ref{HolFactGeneral} taking 
matrix $R$ to be $ C^{\RR} + Q^{\RR} CorrTriDiag$ and $C_1 = C + Q diag(d_n,  ...,d_1)$,
$C_2 = \bar C + \bar Q diag(d_n,  ...,d_1)$, so we get 
${\det}^{\mathrm{col}}_{2n\times 2n}(C^{\RR} + Q^{\RR} CorrTriDiag) = 
{\det}^{\mathrm{col}}_{2n\times 2n}( C + Q diag(d_n,  ...,d_1))  {\det}^{\mathrm{col}}_{2n\times 2n}( \bar C + \bar Q diag(d_n,  ...,d_1))$.

Main theorem \ref{MainTheor} is proved.
\BX

\subsection{Holomorphic factorization for the Capelli matrix \label{SSectHolFactCapel} }
Here we prove the holomorphic factorization theorem for the Capelli matrix, it follows quite easily
from the main factorization theorem from the previous subsection.

{\Th \label{HolFactCapel} Consider elements $E_{ij} = E_{ij}^{Re} + i  E_{ij}^{Im}$, such that they satisfy commutation relations:
$[ E_{ij} , E_{kl} ] = E_{il} \delta_{jk} - E_{kj}\delta_{li} $ i.e. the same commutation relations which satisfy standard basis elements (sometimes called matrix units) in Lie algebra $gl_n$ i.e. matrices with zeros everywhere except $1$ at the position $(ij)$.  Moreover impose relation that all 
$[ E_{ij}, \bar E_{kl} ]=0$, i.e. all elements commute with complex conjugated elements.
Consider matrix: 
\bea
E^{\RR} = \left( \begin{array}{ccccc}
E_{11}^{Re} & E_{11}^{Im} & ... & E_{1n}^{Re} & -E_{1n}^{Im} \\
-E_{11}^{Im} &E_{11}^{Re} & ... & -E_{1n}^{Im} & E_{1n}^{Re} \\
... & ... & ... \\
E_{n1}^{Re} & E_{n1}^{Im} & ... & E_{nn}^{Re} & -E_{nn}^{Im} \\
-E_{n1}^{Im} & E_{n1}^{Re} & ... & -E_{nn}^{Im} & E_{nn}^{Re} \\
\end{array} \right), ~~~
\eea

then
\bea
{\det}^{\mathrm{col}}_{2n\times 2n} (E^{\RR} + CorrTriDiag ) = \\ =
{\det}^{\mathrm{col}}_{n\times n} (E +  diag(n-1, n-2, ...,1,0) ) ~
{\det}^{\mathrm{col}}_{n\times n} (\bar {E}  + diag(n-1, n-2, ...,1,0) ), 
\eea
here $E$ is $n$ by $n$ matrix with $E_{ij}$ on the position $E_{ij}$,
respectively $\bar {E}$ with $\bar E_{ij} = E_{ij}^{Re} - i  E_{ij}^{Im}$;
there are two options for the matrix $CorrTriDiag$ , 
the first one is:
\bea \label{fmlCorrTriDiagCapFactTh}
CorrTriDiag  = \left( \begin{array}{ccccccccc}
n-1 +1/4 & i/4 & 0 & 0 & 0 & 0 & ... & 0 & 0 \\
i/4 & n-1 - 1/4 & 0 & 0 & 0 & 0 & ... & 0 & 0 \\
 0 & 0 & n-2 +1/4 & i/4 & 0 & 0 & ... & 0 & 0 \\
 0 & 0 & i/4 & n-2 - 1/4 &0 & 0 & ... & 0 & 0 \\
 ... & ... & ... & ... & ...  & ... & ... & ... & ... \\
 0 & 0 & 0 & 0 & 0 & 0 & ... & 1/4 & i/4 \\
 0 & 0 & 0 & 0 &0 & 0 & ... & i/4 & -1/4 
\end{array} \right). ~~~
\eea
The second option for $CorrTriDiag$ is the same formula with substitution $i$ to $-i$,
i.e. complex conjugate to the first one.

}

Here ${\det}^{\mathrm{col}}$ is column-determinant.
{See definition \pref{coldetdef}. }

The proof of the theorem immediately follows from the main theorem \ref{MainTheor} and the lemma below.

{\Lem ~}
Denote as usually $\psi^E_k = \sum_i \psi_i E_{ik} $. 
Then it is true that:
\bea
(\psi^E_k)^2 = \psi^E_k \psi_k
\eea
{\bf Proof} of lemma: 
\bea
(\psi^E_k)^2= (\sum_{ij} \psi_i \psi_j E_{ik}E_{jk} ) = 
\sum_{i<j} \psi_i \psi_j [E_{ik},E_{jk}] =\nn\\  = \sum_{j} \psi_k \psi_j [E_{kk},E_{jk}]
= \sum_{j} \psi_k \psi_j (- E_{jk}) 
= 
 (  \sum_{j} \psi_j E_{jk}) \psi_k   =  \psi^E_k \psi_k.
\eea
Lemma is proved. \BX

\section{Real forms  ("decomplexifications ") of the Capelli identities \label{SectCapel} }
This section contains the results on the decomplexification of the Capelli identities.
There are 3 subsections 1) classical identities with square matrices $Z_{ij}$ and $\partial_{ij}$;
and also their symmetric and antisymmetric analogues 2) generalization of \cite{CSS08} identities for 
matrices satisfying certain general commutation relations 3) rectangular case of classical identities.
The strategy of proof is the same - we reduce results to holomorphic factorization identities
from previous section
and known Capelli identities.

\subsection{"Decomplexification" of classical identities}
Here we prove decomplexification results for classical square Capelli identities.
There are two subsections 1) we consider classical case 2) symmetric and antisymmetric analogues.
The proof is the same in all cases.

\subsubsection{Analogues of the Capelli identity \label{SSSectCap} }

Consider the polynomial algebra $\CC[x_{ij},y_{ij}]$, and $2n$ by $2n$  the matrices:
\bea
Z^{\RR} = \left( \begin{array}{ccccc}
x_{11} & y_{11} & ... & x_{1n} & y_{1n} \\
-y_{11} & x_{11} & ... & -y_{1n} & x_{1n} \\
... & ... & ... \\
x_{n1} & y_{n1} & ... & x_{nn} & y_{nn} \\
-y_{n1} & x_{n1} & ... & -y_{nn} & x_{nn} \\
\end{array} \right), ~~~
\DD^{\RR} = (\frac{1}{2}) \left( \begin{array}{ccccc}
\frac{\partial}{\partial x_{11}}  & -\frac{\partial}{\partial y_{11}}  & ...  & \frac{\partial}{\partial x_{n1}} & -\frac{\partial}{\partial y_{n1}} \\
\frac{\partial}{\partial y_{11}}  & \frac{\partial}{\partial x_{11}}  & ...  & \frac{\partial}{\partial y_{n1}} & \frac{\partial}{\partial x_{n1}} \\
... & ... & ... \\
\frac{\partial}{\partial x_{n1}}  & -\frac{\partial}{\partial y_{n1}}  & ...  & \frac{\partial}{\partial x_{nn}} & - \frac{\partial}{\partial y_{nn}} \\
\frac{\partial}{\partial y_{n1}}  & \frac{\partial}{\partial x_{n1}}  & ...  & \frac{\partial}{\partial y_{nn}} & \frac{\partial}{\partial x_{nn}} \\
\end{array} \right). ~~~
\eea

As usually we denote by $D^t$ transpose to matrix $D$.

{\Th \label{ThCap}
\bea
{\det}^{\mathrm{col}} ( Z^{\RR} (\DD^{\RR})^t + \mathrm{CorrTriDiag})= {\det}^{\mathrm{col}} ( Z^{\RR}) {\det}^{\mathrm{col}} ((\DD^{\RR})^t),
\eea
where there are two options for the matrix $CorrTriDiag$ , 
the first one is:
\bea
CorrTriDiag  = \left( \begin{array}{ccccccccc}
n-1 +1/4 & i/4 & 0 & 0 & 0 & 0 & ... & 0 & 0 \\
i/4 & n-1 - 1/4 & 0 & 0 & 0 & 0 & ... & 0 & 0 \\
 0 & 0 & n-2 +1/4 & i/4 & 0 & 0 & ... & 0 & 0 \\
 0 & 0 & i/4 & n-2 - 1/4 &0 & 0 & ... & 0 & 0 \\
 ... & ... & ... & ... & ...  & ... & ... & ... & ... \\
 0 & 0 & 0 & 0 & 0 & 0 & ... & 1/4 & i/4 \\
 0 & 0 & 0 & 0 &0 & 0 & ... & i/4 & -1/4 
\end{array} \right). ~~~
\eea
The second option for $CorrTriDiag$ is the same formula with substitution $i$ to $-i$,
i.e. complex conjugate to the first one. 

}

Here ${\det}^{\mathrm{col}}$ is column-determinant.
{See definition \pref{coldetdef}. }

{\Rem ~} We will work out case $n=1$ in greater generality later (see subsection \pref{n1CapSubSec}).

\PRF
As usually define $z_{ab}= x_{ab} + i y_{ab}$,
 $\partial_{z_{ab}}= (1/2)(\partial_{x_{ab}} - i \partial_{y_{ab}})$;
and matrices $Z,D$ with elements $z_{ij}, \partial_{z_{ij}}$ respectively.
 
Since all elements of $Z$ commute among themselves (and the same for  $\DD$), it is true that:
\bea
{\det}_{2n\times 2n}  ( Z^{\RR} ) = {\det}_{n\times n} (Z){\det}_{n\times n} (\bar Z), \\  
{\det}_{2n\times 2n} ((\DD^{\RR})^t) = {\det}_{n\times n} (\DD^t ) {\det}_{n\times n} (\bar \DD^t ).  
\eea

Hence:
\bea
{\det}_{2n \times 2n}  ( Z^{\RR}) {\det}_{2n \times 2n}((\DD^{\RR})^t) = 
{\det}_{n\times n} (Z) {\det}_{n\times n} (\bar Z) {\det}_{n\times n} (\DD^t) {\det}_{n\times n} (\bar \DD^t)=
\eea 
using that $\bar Z_{ij}$ commute with $ \DD_{kl}$ we get:
\bea
={\det}_{n\times n} (Z)  {\det}_{n\times n} (\DD^t)  {\det}_{n\times n} (\bar Z) {\det}_{n\times n} (\bar \DD^t)
=
\eea 
Using classical Capelli identities:
\bea
={\det}_{n\times n}^{\mathrm{col}} (Z\DD^t +diag(n-1,n-2,...,1,0) )    {\det}_{n\times n}^{\mathrm{col}} (\bar Z\bar \DD^t +diag(n-1,n-2,...,1,0))
=
\eea 
It is well-known that matrix $Z\DD^t$ satisfies the  $gl_n$ commutation relations  (actually its elements are right invariant vector fields on $GL_n$, because $Zgg^{-1}D^t =ZD^t$).
So using our theorem \ref{HolFactCapel} we arrive to the desired result:
\bea
=
{\det}^{\mathrm{col}}_{2n\times 2n} ( Z^{\RR} (\DD^{\RR})^t + \mathrm{CorrTriDiag})
.
\eea
Theorem is proved. \BX

\subsubsection{Decomplexifications of the symmetric and antisymmetric Capelli identities \label{SSSectASymCap}  }
{\Th \label{ThASymCap}
Let $Z,D$ be $n\times n$ matrices which are

either

(Turnbull's case) 
symmetric matrices  defined as: 
 $z_{ij}$ $i\le j$ on both positions $(ij)$ and $(ji)$ and
$\DD$ is the symmetric matrix with  $\partial_{z_{ij}}$ out of the diagonal and $2 \partial_{z_{ii}}$ on the diagonal

or

(Howe-Umeda-Kostant-Sahi case)
$n$ is even and 
$Z$ is the anti-symmetrix  matrix with $z_{ij}$ $i< j$ at the  position $(ij)$ and $-z_{ij}$ at $(ji)$ and respectively $\DD$ is  corresponding anti-symmetric matrix
made of $\partial_{ij}$,

   then the identity holds:  
\bea
{\det}^{\mathrm{col}}_{2n\times 2n} ( Z^{\RR} (\DD^{\RR})^t + \mathrm{CorrTriDiag})= {\det}^{\mathrm{col}}_{2n\times 2n} ( Z^{\RR}) {\det}^{\mathrm{col}}_{2n\times 2n} ((\DD^{\RR})^t).
\eea
}

The proof is the same as in the classical Capelli case.
We should only remark that that in both cases
 elements of the matrix $Z\DD^t$ satisfy
the same $gl_n$ commutatation relations $[E_{ij}, E_{kl}] = \delta_{jk} E_{il} - \delta_{li} E_{kj} $
as it is required in our theorem  \ref{HolFactCapel}.

\subsection{"Decomplexification" of Caracciolo - Sportiello - Sokal identities\label{SSectCCSCap} }

In a recent paper by Caracciolo - Sportiello - Sokal \cite{CSS08}
certain generalizations of  the Capelli identities were proposed.
They have completely algebraic flavor - if certain commutation
relations (let us call them (G,T)-CSS-conditions) between matrix elements of $M,Y,Q$  satisfied , then
$det(M)det(Y) = det(MY+Q~diag(n-1,n-2,...,0))$.
The main result of this section that if we additionally assume that  elements
which belong to the same column of matrix $Y$ i.e. $Y_{*,k}$ 
commute among themselves, then there exists decomplexification of 
CSS-Capelli identities.

Let us give more details.
There are three similar commutation relations which we will 
call CSS, TCSS, GCSS which ensure the existence of the Capelli identities.
Here TCSS (Turnbull-CSS) - for the case of symmetric matrix $M$,
and GCSS (generalized or Grassman) is most general, it subsumes the previous two
as particular cases (it was  proposed in \cite{CFR09}).
We will show that CSS and TCSS combined with  $[Y_{ak}, Y_{bk}]=0$
ensures validity of decomplexified identity.
Unfortunately we are unable to prove the same for GCSS case.
Another unfortunate moment - that we are unable to prove results
for rectangular matrices. 
The reason is that rectangular case involves summation over sub-matrices $C_{IL}$
which are not decomplexification of something, so we do not know how to treat them.
\footnote{In the case of  classical identities we found a way round - we work in commutative case, then
use map "Wick" from commutative algebra to noncommutative  and by some tricks we show that everything works,
however such way does not seem to work in general CSS case since it is not clear what "Wick" can be,
and doubtly it exists at all, since our noncommutative algebra is of much bigger size than commutative one}

Our way of proving Capelli identity is again via the factorization result.
So essentially we need to show that CSS, TCSS combined with  
$[Y_{ak}, Y_{bk}]=0$ ensures that matrices $C=MY$ and $Q$ satisfy 
the conditions in our main theorem. 
About the  GCSS case we can show the validity of the first condition 
in the main theorem. However we do not see the way do deduce the second and the third,
moreover it does not seems plausible for us. It might be that our main theorem is not general
enough to cover GCSS case or it might that GCSS is not enough for factorization and decomplexified
Capelli - all these questions are not clear.

\subsubsection{Manin matrices - reminder \label{SSSectManinMatrReminder}  } 

CSS-Capelli identity naturally involves Manin matrices (\cite{CFR09}).
Manin matrices are certain  matrices with non-necessarily commutative entries.
Main feature about them is that despite non-commutativity, almost all facts of linear algebra
(Cayley-Hamilton, Cramer, Plucker, Schur complement, etc. theorems) can be generalized for them (\cite{CFR09}). 
Any matrix with commutative entries is \MM.
 Let us very briefly remind the definition and one property on \MMs, which will be used.

{\Def \label{D1a22}
Let us call an $n\times n'$ matrix $M$ by {\em \MM},
if the following conditions hold true.
For any $2\times 2$  submatrix
$(M_{ij,kl})$, consisting of rows $i$ and $k$,
and columns $j$ and $l$ (where $1 \leq i < k \leq n$,
and $1 \leq j < l \leq n'$):
\bea
\left(\begin{array}{ccccc}
... & ...& ...&...&...\\
... & M_{ij} &... & M_{il} & ... \\
... & ...& ...&...&...\\
... & M_{kj} &... & M_{kl}& ... \\
... & ...& ...&...&...
\end{array}\right)
\equiv
\left(\begin{array}{ccccc}
... & ...& ...&...&...\\
... & a &... & b& ... \\
... & ...& ...&...&...\\
... & c &... & d& ... \\
... & ...& ...&...&...
\end{array}\right)
\eea
the following commutation relations hold:
\begin{eqnarray}
\label{eq.commutation1}
ac =  ca,~~~ bd =  db, && ~~~ \quad \text{(commutation of the entries in a column)} \\
\label{eq.commutation2}
 ad - da  = cb - bc,&& \quad \text{(cross commutation relation)} .
\end{eqnarray}
}

The definition can be rephrased like this:
\begin{itemize}
\item elements which belong to the same column of $M$ commute among themselves
\item commutators of  cross terms are equal:
$\forall p,q,k,l~~[M_{pq}, M_{kl} ]= [M_{kq}, M_{pl} ]$, \\ e.g. $[M_{11}, M_{22} ]= [M_{21}, M_{12} ]$,
$[M_{11}, M_{2k} ]= [M_{21}, M_{1k} ]$.
\end{itemize}
{\Rem ~} The second condition for the case $q=l$  obviously implies
the first one. 

The following proposition provides both a conceptual point of view on \MMs (\cite{Manin87}),
as well a technical fact which will be often used below.
\MMs~ are "non-commutative endomorphisms" of the polynomial algebra $\CC[x_1,...,x_m]$
(and equivalently of the Koszul dual  $\CC[\psi_1,...,\psi_n]$).
"Non-commutative endomorphisms" of algebra $V$ is a Hopf algebra $A$ such that
exists homomorphism  $\phi : V\to M\otimes V$ (called "coaction") and $\phi$ should satisfy coaction axioms (co-assotiativity like 
relation : $\phi \phi = \delta \phi$). 
(Manin's construction  can be applied to any algebra - what we call Manin matrices
is related to the simplest case: $V=\CC[\psi_1,...,\psi_n]$).

{\Prop \label{Coact-pr} Coaction. Consider a rectangular $n\times
m$-matrix $M$, the polynomial algebra $\CC[x_1,...,x_m]$ and the
Grassmann algebra $\CC[\psi_1,...,\psi_n]$ (i.e. $\psi_i^2=0, \psi_i
\psi_j= -\psi_j \psi_i$); let $x_i$ and $\psi_i$ commute with
$M_{pq}$: $\forall i,p,q:~~[x_i, M_{pq}]=0$, $[\psi_i, M_{pq}]=0$.
Consider new variables $ x_i^M$, $\psi_i^M$: 
\bea
\left(\begin{array}{c}
 x_1^M \\
...  \\
 x_n^M
\end{array}\right)
=
\left(\begin{array}{ccc}
 M_{11} & ... &  M_{1m} \\
... & ... & ... \\
 M_{n1} & ... &  M_{nm}
\end{array}\right)
\left(\begin{array}{c}
 x_1 \\
... \\
 x_m
\end{array}\right),
~~~~
( \psi_1^M, ... ,  \psi_m^M) =
(\psi_1, ... ,   \psi_n)
\left(\begin{array}{ccc}
 M_{11} & ... &  M_{1m} \\
...  & ... & ... \\
 M_{n1} & ... &  M_{nm}
\end{array}\right).
\eea
Then the following three conditions are equivalent:
\begin{itemize}
\item $M$ is a \MM 
\item the variables $ x_i^M$ commute among themselves: $[ x_i^M, x_j^M]=0$
\item the variables $ \psi_i^M$ anticommute among themselves:
$ \psi_i^M \psi_j^M + \psi_j^M \psi_i^M =0$.
\end{itemize}
} {\Rem ~} The conditions $(\psi_i^M)^2=0$ are equivalent to column
commutativity, and  $\psi_i^M \psi_j^M = -\psi_j^M \psi_i^M$, $i<j$, to
the cross term relations.

\subsubsection{(T)CSS-conditions and Capelli identity \label{CCSCapSSSect} }  
Let us give decomplexifications of the Capelli identities found in \cite{CSS08} (we follow notations from \cite{CFR09}).

Consider Grassman algebra $\CC[ \psi_i ]$. As usually for any matrix $N$ we denote by $\psi_k^N= \sum_j \psi_j N_{jk}$,
also $N^\RR$ is decomplexification of $N$ -- i.e. matrix of size $2n\times 2m$, where each element is substituted by 
$real, imag; -imag, real$, algebra $A^\CC$ means scalar extension to "complex numbers" of some algebra $A$ 
(see section \pref{SectPrelim} for further clarifications of notations).

{\Def Let us say that matrices $M,Y,Q$ satisfy CSS-condition if 
\bea
\forall l,j,r,p:~~ [Y_{lj}, M_{rp}] = \delta_{lp} Q_{rj}, \\ 
\mbox{ or equivalently:   } \forall l,j,p:~~
[Y_{lj}, \psi^M_{p}] = \delta_{lp} \psi^Q_{j}.
\eea
}

{\Def Let us say that matrices $M,Y$ satisfy TCSS-condition (Turnbull-CSS) if 
\bea
\forall i,j,k,l:~~ [M_{ij}, Y_{kl}] = -h(\delta_{jk}\delta_{il} + \delta_{ik}\delta_{jl}),
\eea
with element $h$ which is central\footnote{In \cite{CSS08} \cite{CFR09} centrality of $h$ was not required,
here a little strengthen conditions, may be it is better to call this condition as "strict TCSS"}.
}

{\Th \label{CCSCap} Decomplexification of the  (T)CSS  Capelli identities.

Assume  $M,Y,Q$ are square $n\times n$ matrices with elements in $A^\CC$ such that \\
either\\
(1) $M,Y,Q$ satisfy CSS-condition (see above) \\
(2) $M$ is a Manin matrix (see previous subsubsection \pref{SSSectManinMatrReminder})\\
or\\
(1') $M,Y$ satisfy TCSS-conditions (see above), matrix $Q$ is defined as $Q=h~Id$ \\
(2') $M$ is a symmetric matrix with commuting entries\footnote{any symmetric Manin matrix is a matrix with commuting entries if $char(k)\ne 2$, so we may require $M$ is Manin, but it is the same} \\
and in both cases CSS and TCSS we require:\\
(3) elements in the same column of $Y$ commute among themselves i.e. $\forall a,b,k: [Y_{ak}, Y_{bk}]=0$ \\
(4)
for any two elements $\alpha, \beta$ of matrices $M,Y,Q$ it is true that $[\alpha, \bar \beta]=0$, i.e. all elements commute with complex conjugated elements; \\ 
{\bf then}
\bea
{\det}^{\mathrm{col}}_{2n\times 2n} ( M^{\RR} Y^{\RR} + Q^\RR\mathrm{CorrTriDiag})= {\det}^{\mathrm{col}}_{2n\times 2n} ( M^{\RR}) {\det}^{\mathrm{col}}_{2n\times 2n} (Y^{\RR}),
\eea
where
\bea \label{fmlCorrTriDiagInCSS}
CorrTriDiag  = \left( \begin{array}{ccccccccc}
n-1 +1/4 & i/4 & 0 & 0 & 0 & 0 & ... & 0 & 0 \\
i/4 & n-1 - 1/4 & 0 & 0 & 0 & 0 & ... & 0 & 0 \\
 0 & 0 & n-2 +1/4 & i/4 & 0 & 0 & ... & 0 & 0 \\
 0 & 0 & i/4 & n-2 - 1/4 &0 & 0 & ... & 0 & 0 \\
 ... & ... & ... & ... & ...  & ... & ... & ... & ... \\
 0 & 0 & 0 & 0 & 0 & 0 & ... & 1/4 & i/4 \\
 0 & 0 & 0 & 0 &0 & 0 & ... & i/4 & -1/4 
\end{array} \right). ~~~
\eea
}

Here ${\det}^{\mathrm{col}}$ is column-determinant.
{See definition \pref{coldetdef}. }

{\bf Proof} (Modula propositions \pref{propCSSvsFact}, \pref{propTCSSvsFact} ).

The case $n=1$ is somewhat specific and we will consider it separately below. Now we assume $n>1$.

The idea of the proof is the same as for the  other Capelli identities in this paper.
We reduce it to holomorphic factorization main theorem \pref{MainTheor} 
and already known Capelli identities \cite{CSS08}.
So essential part of the proof is to show that (T)CSS conditions combined with
conditions (3),(4) ensure that conditions in the main theorem holds true.
This is given in next subsection  as separated propositions 
\pref{propCSSvsFact}, \pref{propTCSSvsFact} (these propositions require $n>1$). 
These propositions might be of some  interest by themselves.

Let us give a proof modula these propositions.

\bea
{\det}^{\mathrm{col}}_{2n\times 2n} ( M^{\RR}) {\det}^{\mathrm{col}}_{2n\times 2n} (Y^{\RR})
=
\eea
Since elements in the same column of matrices $M,Y$ commute among themselves,
and also real parts of any element commute with imaginary parts of any element,
we can apply our weak factorization theorem \pref{theor1} or our main theorem \pref{MainTheor}
and write the factorizations:
\bea
=
{\det}^{\mathrm{col}}_{n\times n} ( M) {\det}^{\mathrm{col}}_{n\times n} ( \widebar{ M} )
{\det}^{\mathrm{col}}_{n\times n} ( Y) {\det}^{\mathrm{col}}_{n\times n} (  \widebar{ Y} )
 =
\eea
Use that "bar" elements commute with non-bar:
\bea
=
{\det}^{\mathrm{col}}_{n\times n} ( M) {\det}^{\mathrm{col}}_{n\times n} ( Y)
 {\det}^{\mathrm{col}}_{n\times n} ( \widebar{ M} )  {\det}^{\mathrm{col}}_{n\times n} (  \widebar{ Y} )
 =
\eea
Use Capelli identities from \cite{CSS08}:
\bea
=
{\det}^{\mathrm{col}}_{n\times n} ( MY + Q~diag(n-1,n-2,...,0))
 {\det}^{\mathrm{col}}_{n\times n} ( \widebar{ M  Y} + \widebar{Q}~diag(n-1,n-2,...,0)  )
 =
\eea
Use our main theorem \pref{MainTheor} (here we need propositions 
\pref{propCSSvsFact}, \pref{propTCSSvsFact} (these propositions require $n>1$))
described below which ensure
that conditions in our main theorem (reminded below \ref{psiCsq}-\pref{psiq2eq0})
follows from (T)CSS conditions plus conditions (3),(4)):

\bea
=
{\det}^{\mathrm{col}}_{2n\times 2n} ( M^\RR Y^\RR + Q^\RR~CorrTriDiag).
\eea

Proof finished \BX.

\subsubsection{Example $n=1$ \label{n1CapSubSec} } 

The case $n=1$ is specific and does not demonstrate main features of the general case 
(the correction matrix $Q$ may not be related to commutators of $M,Y$, it can be zero; 
also the key propositions \pref{propCSSvsFact}, \pref{propTCSSvsFact}
use  require $n>1$).
Nevertheless it might be worth to work out $n=1$ case explicitly.

{\Prop Consider three elements $M,Y,Q$ such that any of them commute with  complex conjugation of another
(i.e.  $[\alpha, \bar \beta ] =0$, for $\alpha, \beta$ is any of $M,Y,Q$), then
\bea
{\det}^{\mathrm{col}}_{2\times 2}
\left( \left( \begin{array}{cc}
Re(M)  & Im(M) \\
-Im(M)  & Re(M) \\
\end{array} \right)
\left( \begin{array}{cc}
Re(Y)  & Im(Y) \\
-Im(Y)  & Re(Y) \\
\end{array} \right)
+
\left( \begin{array}{cc}
Re(Q)  & Im(Q) \\
-Im(Q)  & Re(Q) \\
\end{array} \right)
\left( \begin{array}{cc}
1/4  & i/4 \\
i/4  & -1/4 \\
\end{array} \right)
\right) = \\ =
{\det}^{\mathrm{col}}_{2\times 2}
 \left( \begin{array}{cc}
Re(M)  & Im(M) \\
-Im(M)  & Re(M) \\
\end{array} \right)
{\det}^{\mathrm{col}}_{2\times 2}
\left( \begin{array}{cc}
Re(Y)  & Im(Y) \\
-Im(Y)  & Re(Y) \\
\end{array} \right).
\eea 
}

{\Rem ~} Comparing to the general $n$  case, here we see that  $Q$ is arbitrary and not related to the commutators of $M,Y$.
Pay attention that the left hand side contains $Q$ while right hand side does not,
this is somewhat surprising, but indeed true. 

\PRF

Using proposition \pref{n1pr2} we can transform the left hand side to the form:

\bea
{\det}^{\mathrm{col}}_{2\times 2}
\left( \left( \begin{array}{cc}
Re(M)  & Im(M) \\
-Im(M)  & Re(M) \\
\end{array} \right)
\left( \begin{array}{cc}
Re(Y)  & Im(Y) \\
-Im(Y)  & Re(Y) \\
\end{array} \right)
\right)=
\eea
This can be transformed to:
\bea
={\det}^{\mathrm{col}}_{2\times 2}
 \left( \begin{array}{cc}
Re(MY)  & Im(MY) \\
-Im(MY)  & Re(MY) \\
\end{array} \right) = Re(MY)^2 + Im(MY)^2 =
\eea
Using that $[MY, \widebar{MY}] = 0$, we get:
\bea
= MY\widebar{MY}  =
\eea
Using that $[Y, \widebar{M}] = 0$, and $[M, \bar M] = [Y, \bar Y]=0$ we get:
\bea
 =M\bar M Y\widebar{Y}  = (Re(M)^2+ Im(M)^2)( Re(Y)^2+ Im(Y)^2)  = \\ =
{\det}^{\mathrm{col}}_{2\times 2}
 \left( \begin{array}{cc}
Re(M)  & Im(M) \\
-Im(M)  & Re(M) \\
\end{array} \right)
{\det}^{\mathrm{col}}_{2\times 2}
\left( \begin{array}{cc}
Re(Y)  & Im(Y) \\
-Im(Y)  & Re(Y) \\
\end{array} \right).
\eea 
The proof is finished.
\BX

\subsubsection{(T)CSS-conditions vs factorization conditions} 
Let us show that (T)CSS-conditions and condition (3) (i.e. commutativity 
of elements $Y_{*k}$ from the same column $k$ of matrix $Y$) ensures the conditions which
we used in our main theorem.   These propositions have been used in the proof of the theorem above,
and also can be of interest by themselves.

Let us first remind the conditions on $M,Y,Q$ which we have in our main theorem \pref{MainTheor}
:

{\bf "Factorization conditions" for matrices $C$ and $Y$} (see main theorem \pref{MainTheor}):
\bea \label{psiCsq}
\forall k~~~~ (\psi^C_k)^2 =  \psi^C_k \psi^Q_k  ,
\mbox{ or equivalently:   }
\forall i,j,k: ~~[C_{ik}, C_{jk}] = C_{ik} Q_{jk} -  C_{jk}Q_{ik},
 \\
\label{psiCpsiQ}
\forall k~~~~ \psi^C_k \psi^Q_k + \psi^Q_k \psi^C_k = 0 , \mbox{ or equivalently:   }
[C_{ik}, Q_{jk}] = [C_{jk},  Q_{ik}], \\
 \label{psiq2eq0}
\forall k~~~~ (\psi^Q_k)^2 = 0, \mbox{ or equivalently:   }
 [Q_{ik}, Q_{jk}] = 0.
\eea

{\Prop \label{propCSSvsFact} Assume that:
\begin{enumerate}
\item matrix $M$ 	is a Manin matrix;
\item matrices $M,Y$ satisfy the CSS-condition 
$[Y_{lj}, \psi^M_{p}] = \delta_{lp} \psi^Q_{j}$;
\item matrix $Y$ is such that elements
from the same column commute among themselves (i.e. $\forall a,b,k~~[Y_{ak}, Y_{bk}]=0$ ;
\item char (k) not equal to 2, the number of columns in $M$ is greater than 1; 
\end{enumerate}

{\bf then} \\
 the "factorization" conditions \ref{psiCsq}-\ref{psiq2eq0} for matrices $C=MY$ and $Q$ required in main theorem \pref{MainTheor} will be satisfied.
}

Before proof of this proposition let us similar one related to the case of TCSS condition.

{\Prop \label{propTCSSvsFact} Assume $M$ is symmetric matrix with commuting entries,
and matrices $M$ and $Y$ satisfy the TCSS (Turnbull-CSS) condition:
\bea
[M_{ij}, Y_{kl}] = -h(\delta_{jk}\delta_{il} + \delta_{ik}\delta_{jl}),
\eea
with element $h$ which is central,

{\bf then}\\
 all the conditions \ref{psiCsq}-\ref{psiq2eq0} for matrices $C=MY$ and $Q$ required in main theorem \pref{MainTheor} will be satisfied.
}

\PRF Let us prove proposition \ref{propCSSvsFact} above.

Let us prove formula \ref{psiCpsiQ} : $\psi^C_j \psi^Q_j + \psi^Q_j \psi^C_j = 0$.

\bea
\psi^C_j \psi^Q_j = 
\eea
Obviously $\psi^C_j = \sum_i \psi^M_i Y_{ij}$.
CSS-condition implies:
$\psi^Q_j = [Y_{aj},\psi^M_{a} ]$,
so we can write:
\bea
= (\sum_i \psi^M_i Y_{ij} ) [Y_{aj},\psi^M_{a} ] =
\eea
Let us choose $a\ne i$ (we can do this since number of columns in $M$ is bigger than 1 - it was required). 
We see that 
\bea
= - [Y_{aj},\psi^M_{a} ] (\sum_i \psi^M_i Y_{ij} )  = \psi^Q_j \psi^C_j .
\eea
 because elements super-commute: 
$[Y_{aj}, Y_{ij} ] =0$ (requirement in theorem), 
$0 = [Y_{ij},\psi^M_{a} ] = [Y_{aj},\psi^M_{i} ]$ (from CSS-condition),
$\psi^M_{a} \psi^M_i = - \psi^M_i \psi^M_{a} $ (from Manin's property).

Formula \ref{psiCpsiQ} is proved. \\

Let us prove formula \ref{psiq2eq0} : $(\psi_k^Q)^2=0$. 

Proof uses the same idea.
Take $l\ne p$ (we can do this since number of columns in $M$ is bigger than 1 - it was required), CSS-condition implies 
$\psi^Q_j = [Y_{lj},\psi^M_{l} ] = [Y_{pj},\psi^M_{p} ] $
and $0 = [Y_{pj},\psi^M_{l} ] = [Y_{lj},\psi^M_{p} ]$.
From the second one, and from $[Y_{lj}, Y_{pj} ] =0$, and from Manin's property 
$\psi^M_{l} \psi^M_p = - \psi^M_p \psi^M_{l} $
we can conclude 
$[Y_{lj},\psi^M_{l} ]  [Y_{pj},\psi^M_{p} ] = - [Y_{pj},\psi^M_{p} ][Y_{lj},\psi^M_{l} ]  $
since all elements are super-commuting.

On the other hand 
$(\psi^Q_j)^2  = [Y_{lj},\psi^M_{l} ]  [Y_{pj},\psi^M_{p} ] = - [Y_{pj},\psi^M_{p} ][Y_{lj},\psi^M_{l} ] = - (\psi^Q_j)^2  $.
So we conclude that $(\psi^Q_j)^2  =0$.

Formula \ref{psiq2eq0} is proved.

Let us prove formula \ref{psiCsq} : $(\psi^C_j)^2  =  \psi^C_j \psi^Q_j $.\footnote{We will prove later
more general fact, that GCSS condition implies \ref{psiCsq} (proposition \pref{PropGCSSfact}).
It was proved in \cite{CFR09} that CSS implies GCSS.} 
\bea
(\psi^C_j)^2  = \sum_i \psi^M_i Y_{ij} \sum_k \psi^M_k Y_{kj}=
\eea
Permute $Y_{ij}$ and $\psi^M_k$:
\bea
=\sum_i\sum_k \psi^M_i \psi^M_k Y_{ij}    Y_{kj} +
\sum_i\sum_k \psi^M_i [Y_{ij},  \psi^M_k]  Y_{kj} =
\eea
The first sum is zero because $\psi^M_i, \psi^M_k$ anticommute (Manin's property), while $Y_{ij},    Y_{kj}$ commute (requirement in theorem), so it follows from the $Tr(AB)=0$ for  $A$ is antisymmetric, while $B$ antisymmetric, which can be seen like this: $\sum_{ij} A_{ij}B_{ji} = $ use symmetry $= \sum_{ij} -A_{ji}B_{ij}=$ denote $i$ by $j$ and vice versa $=  \sum_{ji} -A_{ij}B_{ji}$, so we get same expression with minus sign, so $\sum_{ij} A_{ij}B_{ji} = 0$, if $char \ne 2$.  
For the second sum we use CSS-condition:
\bea
= 0+
\sum_i\sum_k \psi^M_i \delta_{ik} \psi^Q_j   Y_{kj}
= \sum_k \psi^M_k  \psi^Q_j   Y_{kj} = 
\eea
It was proved in \cite{CFR09} (section 6.3 lemma 20 page 44) that CSS-condition implies:  $\psi^M_k  \psi^Q_j =  - \psi^Q_j \psi^M_k$, so:
\bea
= - \psi^Q_j \sum_k  \psi^M_k     Y_{kj} = - \psi^Q_j  \psi^C_j =
\eea
finally by previously proved formula 
\ref{psiCpsiQ} : $\psi^C_j \psi^Q_j + \psi^Q_j \psi^C_j = 0$, so
we get:
\bea
=   \psi^C_j \psi^Q_j. 
\eea
Formula \ref{psiCsq} is proved.

Proposition is proved. \BX

\PRF Let us prove proposition \ref{propTCSSvsFact} above.

Conditions \ref{psiCpsiQ}, 
\ref{psiq2eq0}
are satisfied by trivial reason - $h$ is central.
Instead of giving the direct proof of the condition  \ref{psiCsq}
let us refer to the more general proposition  \pref{PropGCSSfact}
which we will prove in the next subsection.
It states that GCSS condition implies desired  condition \ref{psiCsq},
while 
in the \cite{CFR09} section 6.4 page 47 it was proved that TCSS conditions
imply GCSS conditions. h

Proposition is proved.\BX

\subsubsection{GCSS-conditions } 
In \cite{CFR09} there was proposed a condition which generalizes both CSS and TCSS conditions,
let us call it GCSS (Grassman or generalized).
Here we remind it, relate it to our factorization conditions and give some reformulation.

{\Def Let us say that matrices $M,Y,Q$ satisfy GCSS-condition (Grassman or generalized CSS) if 
\bea
(1) ~~~~~ \sum_l \psi^M_l[Y_{lj}, \psi^M_p] =  \psi^M_p  \psi^Q_j, \\
(2) ~~~~~ \psi^Q_j\psi^M_p + \psi^M_p \psi^Q_j = 0.
\eea
}

{\Prop Assume that: \label{PropGCSSfact}
\begin{enumerate}
\item matrix $M$ 	is a Manin matrix;
\item matrices $M,Y$ satisfy the GCSS-condition above 
\item matrix $Y$ is such that elements
from the same column commute among themselves (i.e. $\forall a,b,k~~[Y_{ak}, Y_{bk}]=0$ ;
\item char (k) not equal to 2, the number of columns in $M$ is greater than 1; 
\end{enumerate}

{\bf then} \\
 the first "factorization" condition \ref{psiCsq} is true for matrices $C=MY$ and $Q$,
 i.e. 
 \bea
 \forall k~ (\psi^C_k)^2 =  \psi^C_k \psi^Q_k.
 \eea

}

\PRF
Let us multiply the equality \ref{GCSS1} by $Y_{pj}$ and make summation over $p$: $0=\sum_p \psi^C_j\psi^M_p Y_{pj} + \sum_p  \psi^M_p (\psi^C_j - \psi^Q_j)Y_{pj} = 
\psi^C_j\psi^C_j +   \sum_p  \psi^M_p \psi^C_jY_{pj}  + \sum_p  \psi^M_p\psi^Q_j Y_{pj} $ by the second condition \ref{GCSS2}
we can permute $\psi^M$ and $\psi^Q$:

$ \psi^C_j\psi^C_j +   \sum_p \sum_a \psi^M_p \psi^M_a Y_{aj} Y_{pj}  + \sum_p  \psi^Q_j\psi^M_p Y_{pj}= \psi^C_j\psi^C_j +   0  +  \psi^Q_j\psi^C_j $.

Where we used $\sum_p \sum_a \psi^M_p \psi^M_a Y_{aj} Y_{pj}=0$,
since $\psi^M$ are anticommuting by Manin's property and $Y_{*j}$
are commuting by requirement.
Comparing initial zero at the left hand side and the last expression we get the desired result.

Proposition is proved. \BX

{\Prop Assume $M$ is a \MM~. GCSS conditions: 
can be reformulated as follows:
\bea
(1) ~~~~~ \psi^C_j\psi^M_p + \psi^M_p (\psi^C_j - \psi^Q_j) = 0, \label{GCSS1}\\
(2) ~~~~~ \psi^Q_j\psi^M_p + \psi^M_p \psi^Q_j = 0.\label{GCSS2}
\eea
}

\PRF Indeed, by requirement: $\psi^M_p  \psi^Q_j = \sum_l \psi^M_l[Y_{lj}, \psi^M_p]$
let us transform right hand side $=
\sum_l \psi^M_l Y_{lj} \psi^M_p - \sum_l \psi^M_l \psi^M_p Y_{lj}  =
 \psi^C_j  \psi^M_p + \sum_l \psi^M_p \psi^M_j  Y_{lj}  =
 \psi^C_j  \psi^M_p + \psi^M_p \psi^C_j $.
Comparing the first one and the last one we get the desired equality \ref{GCSS1}.

Proposition proved. \BX

\subsection{Rectangular matrices}
For  matrices with commutative elements there is the well-known
Cauchy-Binet formula: $\det((AB)_{IJ})=\sum_{L} \det(A_{IL})\det (B_{LJ})$,
where $I,L,J$ are multi-indexes. There are generalizations of the Capelli identity
in this direction.
Here we prove "decomplexifications"
of the Capelli-Cauchy-Binet formulas for standard Capelli matricies,
as well as for symmetric ones, we also formulate conditional result in
the antisymmetric case.\footnote{In antisymmetric case we are not sure that the
Capelli-Cauchy-Binet formula is known without decomplexification. If it is known
then our argument proves its decomplexified analog. So result is conditional.} 
The subtlety is that submatrices $A_{IL},B_{LJ}$ are not decomplexifications 
of some submatrices this makes the proof not so simple as for square matrices.
This subtlety is the reason why we are not able to decomplexify the rectangular
case of the CSS-type identities discussed in the previous section.

\subsubsection{Analogues of Capelli's and  Turnbull's (symmetric) cases \label{SSSectRectCapTur} }
 As usually for any matrix $N$ we denote by $N^\RR$ its decomplexification of -- i.e. matrix of size $2n\times 2m$, where each element is substituted by 2x2 submatrix
$real, imag; -imag, real$
(see section \pref{SectPrelim} for further clarifications of notations).
We also denote by $D^t$ transpose to matrix $D$.

{\Th \label{ThRectCapTur}
Let $Z,D$ be 
either

(Capelli's case) 
$Z$ is a matrix with $z_{ij}$ at place  $(ij)$, and $\DD$ with $\partial_{z_{ij}}$ 

or 

(Turnbull's case) $Z,D$ are  
symmetric matrices  defined as: 
 $z_{ij}$ $i\le j$ on both positions $(ij)$ and $(ji)$ and 
$\DD$ is the symmetric matrix with  $\partial_{z_{ij}}$ out of the diagonal and $2 \partial_{z_{ii}}$ on the diagonal

Consider multi-indexes $I=(i_1,...,i_{r})$,
$J=(j_1,...,j_{r})$. The following is true:

\bea
{\det}_{2r\times 2r} ( (Z ~\DD^t)^{\RR}_{IJ} +Q^{\RR} \mathrm{CorrTriDiag})= 
\sum_{L=(l_1<l_2<l_3<...<l_{2r})}
{\det}_{2r\times 2r} ( Z^{\RR}_{IL})
 {\det}_{2r\times 2r} ((\DD^{\RR})^t_{LJ}).
\eea
Where matrices $Q$ and $CorrTriDiag$ are defined
\bea
(Q)_{ab} = \delta_{i_a j_b} ,
\eea

\bea
CorrTriDiag  = \left( \begin{array}{ccccccccc}
r-1 +1/4 & i/4 & 0 & 0 & 0 & 0 & ... & 0 & 0 \\
i/4 & r-1 - 1/4 & 0 & 0 & 0 & 0 & ... & 0 & 0 \\	
 0 & 0 & r-2 +1/4 & i/4 & 0 & 0 & ... & 0 & 0 \\
 0 & 0 & i/4 & r-2 - 1/4 &0 & 0 & ... & 0 & 0 \\
 ... & ... & ... & ... & ...  & ... & ... & ... & ... \\
 0 & 0 & 0 & 0 & 0 & 0 & ... & 1/4 & i/4 \\
 0 & 0 & 0 & 0 &0 & 0 & ... & i/4 & -1/4 
\end{array} \right). ~~~
\eea
The second option for $CorrTriDiag$ is the same formula with substitution $i$ to $-i$.
}

Here ${\det}^{\mathrm{col}}$ is column-determinant.
{See definition \pref{coldetdef}. }

{\Rem ~} Let us emphasize that the theorem above concerns not all submatrices of $(ZD^t)^\RR$,
but only "complex submarices" i.e. ones which have a form $M^\RR$ for some $M\subset (ZD^t)$.
We are not aware whether the theorem might be
extended to all submatrices of $(ZD^t)^\RR$, not just "complex ones"  (in commutative case, it is, of course, so).
Let us mention, that the summation in multi-index $L$ runs through all possible submatrices i.e. 
$Z^{\RR}_{IL}$ $(\DD^{\RR})^t_{LJ}$ do not have a form $M^\RR$. 

\PRF
Let us consider commutative variables $p^x_{ij}, p^y_{ij}$ which will play a role of symbols of  	 
$\frac{\partial}{\partial x_{ij}}, \frac{\partial}{\partial y_{ij}}$.
Denote by $P,P^\RR$ respectively  matrices analogous to $\DD, \DD^\RR$.

Denote by "Wick" the map from the commutative polynomial algebra generated by $ x_{ij}, y_{ij}, p^x_{ij}, p^y_{ij}$ to non-commutative algebra of differential operators generated by
$ x_{ij}, y_{ij}, \frac{\partial}{\partial y_{ij}} , \frac{\partial}{\partial x_{ij}}$  
defined by the rule $x,y$ stays first, derivatives go second, i.e. linear map defined on monomials as follows: 
\bea
Wick: ~~ \prod_{ij} x_{ij}^{n_{x,ij}} \prod_{ij} y_{ij}^{n_{y,ij}} 
\prod_{ij} (p_{ij}^x)^{n_{p_x, ij}} \prod_{ij} (p_{ij}^y)^{m_{p_y, ij}} 
~\mapsto ~
\prod_{ij} x_{ij}^{n_{x,ij}} \prod_{ij} y_{ij}^{n_{y,ij}} 
\prod_{ij} ( \frac{\partial}{\partial x_{ij}} )^{n_{p_x, ij}} 
\prod_{ij} (\frac{\partial}{\partial y_{ij}} )^{m_{p_y, ij}} .
\eea

The usual commutative Cauchy-Binet formula gives:
\bea \label{comCBloc111}
\sum_{L=(l_1<l_2<l_3<...<l_{2r})}
{\det}_{2r\times 2r} ( Z^{\RR}_{IL})
 {\det}_{2r\times 2r} ((P ^{\RR})^t_{LJ}) = 
{\det}_{2r\times 2r} ( (Z ~P^t)^{\RR}_{IJ} ) = 
\eea

Using commutative holomorphic factorization  we can  continue
\bea
= {\det}_{r\times r} ( (Z ~P^t)_{IJ} ) 
{\det}_{r\times r} ( (\bar Z \bar P^t)_{IJ} ) =
\eea

Use commutative Cauchy-Binet twice:

\bea \label{aftercomCB}
=
\sum_{L=(l_1<l_2<l_3<...<l_{r})}
{\det}_{r\times r} ( (Z )_{IL} ) {\det}_{r\times r} ( (P^t)_{LJ} )
\sum_{K=(k_1<k_2<k_3<...<k_{r})}
{\det}_{r\times r} ( (\bar Z )_{IK} ) {\det}_{r\times r} ( (\bar P^t)_{KJ} )
\eea

So comparing \ref{comCBloc111} and \ref{aftercomCB}  we get:
\bea
\sum_{L=(l_1<l_2<l_3<...<l_{2r})}
{\det}_{2r\times 2r} ( Z^{\RR}_{IL})
 {\det}_{2r\times 2r} ((P^{\RR})^t_{LJ}) = \\ =
\sum_{L=(l_1<l_2<l_3<...<l_{r})}
{\det}_{r\times r} ( (Z )_{IL} ) {\det}_{r\times r} ( (P^t)_{LJ} )
\sum_{K=(k_1<k_2<k_3<...<k_{r})}
{\det}_{r\times r} ( (\bar Z )_{IK} ) {\det}_{r\times r} ( (\bar P^t)_{KJ} ) 
\eea

Applying  Wick map to the equality above we get:

\bea \label{fmlLoc1111}
\sum_{L=(l_1<l_2<l_3<...<l_{2r})}
{\det}_{2r\times 2r} ( Z^{\RR}_{IL})
 {\det}_{2r\times 2r} ((\DD^{\RR})^t_{LJ}) 
 = \\ =
\sum_{L=(l_1<l_2<l_3<...<l_{r})}
\sum_{K=(k_1<k_2<k_3<...<k_{r})}
{\det}_{r\times r} ( (Z )_{IL} ) {\det}_{r\times r} ( (\bar Z )_{IK} )
 {\det}_{r\times r} ( (\DD^t)_{LJ} )
 {\det}_{r\times r} ( (\bar \DD^t)_{KJ} ) =
\eea
using commutativity between $\bar Z$ and $\DD$ we rewrite:
\bea
=\sum_{L=(l_1<l_2<l_3<...<l_{r})}
{\det}_{r\times r} ( (Z )_{IL} )  {\det}_{r\times r} ( (\DD^t)_{LJ} )
\sum_{K=(k_1<k_2<k_3<...<k_{r})}
 {\det}_{r\times r} ( (\bar Z )_{IK} )
 {\det}_{r\times r} ( (\bar \DD^t)_{KJ} ) =
\eea

We can apply Capelli-Cauchy-Binet formula (in symmetric case such formula seems first been proven in 
\cite{CSS08} (see also \cite{CFR09}) extending Turnbull's result to rectangular matrices, we are not aware of any proof in antysmmetric case): 
\bea
=
{\det}_{r\times r} ( (Z\DD^t )_{IJ} + Q diag(r-1,r-2,...,1,0)  )
 {\det}_{r\times r} ( (\bar Z \bar \DD^t )_{IJ} +  Q diag(r-1,r-2,...,1,0))
\eea

The matrix $(\bar Z \bar \DD^t )_{IJ}$ is a submatrix of Capelli type
matrix, it is quite easy to see
\footnote{the same arguments as in derivation of theorem
\pref{HolFactCapel} about the holomorphic factorization for the  Capelli
matrix from the main theorem  \pref{MainTheor}. With
the only difference that we deal with submatrix of Capelli matrix, not
the whole matrix.}
 that we can apply our main 
theorem \ref{MainTheor}
and we get
\bea \label{fmlLoc2222}
=
{\det}_{2r\times 2r} ( (Z ~\DD^t)^{\RR}_{IJ} +Q^{\RR} \mathrm{CorrTriDiag}).
\eea
Comparing \ref{fmlLoc1111} and  \ref{fmlLoc2222}
we see the desired equality. 
Theorem is proved.
\BX

\subsubsection{Conditional result in antisymmetric case \label{SSSectRectAntiSym} }
We are not aware of the proof of the Capelli-Cauchy-Binet formula
for rectangular matrices in antisymmetric case.
That it is why we are not able to extend our previous result to this case.
So we formulate the conditional result.

{\bf Conditional Theorem} Consider the setup of the previous theorem,
assuming that $Z$ and $D$ are antisymmetric matrices with elements
$z_{ij}$ and $\partial_{ij}$ on corresponding positions. 
Then the Capelli-Cauchy-Binet formula (the same as in the previous theorem) holds true for $Z^{\RR}$ and $D^{\RR}$,
conditionally that Capelli-Cauchy-Binet formula holds true for submatrices of standard (non-decomplexified) matrices $Z$ and $D$.

Proof is the same as in the previous section.

{\Rem ~ \label{HUKSisHardRemark} } Let us mention that we tried to obtain the proof of the Capelli 
result for antisymmetric matrices (i.e. reprove Howe-Umeda Kostant-Sahi result),
by the methods of \cite{CFR09}, but we failed. Papers \cite{FZ93} (page 9 bottom), \cite{CSS08} (page 9 top) contain similar remarks.
Moreover in \cite{CSS08} conjectures were proposed
(page 36 bottom conjectures 5.1, 5.2)
 extending the
result to $(2n-1)\times (2n-1)$ antisymmetric matrices. They are still open to the best of our knowledge. So there is certain difficulty to obtain combinatorial or Grassman algebra proof of the identity for antisymmetric matrices. We think that it is related to some non-local character of the cancellation of unwanted terms. It would be desirable to find
some simple direct proof of this identity, since it may lead 
to understanding how such non-local cancellations may happen, which hopefully
be helpful for further generalizations.

\section{Cayley identities \label{SectCayl} } 
The so-called Cayley identity $det(d/dx_{ij} )det(x_{ij} )^s = s(s+1)...(s+n-1) det(x_{ij} )^{s-1}$ is a beautiful fact from nineteenth century, alive nowadays. We refer to remarkable paper \cite{CSS11} for all (various proofs of old and new results, 
extensive reference list, history account, and relations to the Capelli identity).

Here we remark that "decomplexification" and "dequaternification"
can be applied to Cayley identities and it is quite direct. Before doing this let us mention two facts
on wider context.

{\bf Cayley from Capelli.} Cayley identity can be easily  derived from the Capelli identity and certain facts on representation theory. In particular let us stress that the correction matrix $diag(n-1,n-2,...,1,0)$
in the Capelli identity exactly corresponds to multiplier $s(s+1)(s+2) ... (s+n-1)$ in the Cayley identity.
The derivation can be sketched like this: $det(D)det(X)^s= det(D^tX+diag(n-1,n-2,...,1,0)) det(X)^{s-1}=$ - we used Capelli identity; now we can use that all out-diagonal terms of $(D^tX+diag(n-1,n-2,...,1,0))$
act by zero on the $det(X)^{s-1}$ because all the subgroup $SL(n)\in GL(n)$ acts trivially on it
so we get: $=det(diag( D^tX+diag(n-1,n-2,...,1,0))) det(X)^{s-1}=$;
now diagonal terms acts as $GL(1)\in GL(n)$ subgroup of scalar matrices modula shifts by $n-k+1$,
so we finally get $s(s+1)...(s+n-1) det(x_{ij} )^{s-1}$.

Discussing the Cayley identity it is worth to mention that it provides 

{\bf Cayley and Bernstein-Sato.} Discussing the Cayley identity it is worth to mention that it provides explicit and simple example for what is called Bernstein polynomial.
Consider any polynomial $P(x_1,...,x_n)$ in $n$-variables.
Berntein proved \cite{B72} that there exists  a polynomial-coefficient partial
differential operator $Q(s, x_i, \partial_{x_i})$ and a polynomial $b(s)\ne 0$   satisfying
\bea
Q(s,x_i, \partial_{x_i}) P(x)^s = b(s) P(x)^{s-1}. 
\eea

In general it is quite difficult to find Bernstein polynomial $b(s)$ (also called Bernstein-Sato polynomial). The roots of it are algebraic numbers related to monodromy of the singularity given by $P(x)=0$.
The Cayley identity provides non-trivial example where $b(s)$  is given explicitly.

\subsection{"Decomplexification" \label{SectCaylDecomp} } 

{\Th Consider $n\times n$ matrices $Z$ and $D$, with elements $z_{ij}$ and $\partial_{z_{ij}}$, let $D^\RR, Z^\RR$ be their "decomplexifications", 
then
\bea
\det (D^\RR) (\det(Z^\RR) )^s = b(s)^2  (\det(Z^\RR) )^{s-1}, 
\eea
where $b(s) =s(s+1)...(s+n-1)$.
}

\PRF
Elements of $Z^\RR$  commute among themselves,
same for $D^\RR$. So $\det(Z^\RR) = \det(Z) \det(\bar Z)$, same for $D^\RR$.
Hence 
\bea
\det (D^\RR) (\det(Z^\RR) )^s =
 \det (D)  \det (\bar D)  (\det(Z) )^s  (\det(\bar Z)  )^s
=  
\eea
Since $\det (\bar D)$ commute with $(\det(Z) )^s$
we can continue
\bea
= \det (D)   (\det(Z) )^s \det (\bar D)   (\det(\bar Z)  )^s
=b(s) (\det(Z) )^{s-1} b(s) (\det(\bar Z) )^{s-1}
= b(s)^2 (\det(Z^\RR) )^{s-1}.
\eea
Theorem is proved. 
\BX

{\Th  Similar "decomplexification" result holds true for symmetric and antisymmetric
complex matrices. 
}

We refer to \cite{CSS11} for symmetric and antisymmetric versions of the 
Cayley identity. The proofs goes exactly in the same way as above.

{\Rem ~} There are more general forms of the Cayley identities. Instead of
$\det(D)$ one considers $\det(Z_{IJ})$ (\cite{CSS11} formula 2.2 page 7,
2.4 page 7 (in symmetric case), etc.)
Clearly the "decomplexification" results extend to this case.

\subsection{"Dequaternionification" \label{SectCaylDecomp} } 

Quaternions can be seen as $\CC^2$ or as $\RR^4$.
So a quaternionic matrix $M$ has its real form $M^\RR$ and 
its complex form $M^\CC$.
The complex form is made of substituting each quaternion $ q=x_0+ix_1+jx_2+kx_3=z_1 + jz_2$ by the 2x2 matrix
\bea
\left( \begin{array}{ccccccccc}
z_1 & z_2\\
-\bar z_2 &\bar z_1\end{array} \right).
\eea
Respectively $M^\RR$ is decomplexification of this matrix.
It is well-known that $\det(M^\RR) = \det(M^\CC)^2$.
For example for 1x1 matrix $M=q=x_0+ix_1+jx_2+kx_3=z_1 + jz_2$,
we have $det(M^\CC) = z_1\bar z_1 + z_2\bar z_2 = x_0^2+x_1^2+x_3^2+x_4^2$.
While $det(M^\RR) = (z_1\bar z_1 + z_2\bar z_2)^2 = (x_0^2+x_1^2+x_3^2+x_4^2)^2$.

{\bf Observation} {\em Consider  $n\times n$ quaternionic matrices $Z$, $D$, with elements 
$q_{ab} = z_{1,ab}+jz_{2,ab}$ and $d_{ab} = 1/2(\partial_{z_{1,ab}}-j\partial_{z_{2,ab}})$ 
and $Z^\CC$, $D^\CC$ their complex forms.

Then commutation relations between elements 
$Z^\CC_{uv}$ $u,v=1...2n$  and $2(D^\CC_{tr})^t$, $r,s=1...2n$ 
are canonical ones:
\bea
[Z^\CC_{uv}, 2 D^\CC_{rs}] = \delta_{ur} \delta_{vs}.
\eea
}

Hence the standard Cayley identity can be applied to matrices $Z^\CC$, $2(D^\CC)^t$
and so we get:

{\Cor  "Dequaternionified" Cayley 1.
\bea
\det(D^\CC) \det(Z^\CC)^s  =\frac{1}{2^{2n}} \det(2(D^\CC)^t) \det(Z^\CC)^s  =  \frac{1}{2^{2n}} s(s+1)...(s+2n-1)\det(Z^\CC)^{s-1}
\eea
}

{\Cor  "Dequaternionified" Cayley 2.
\bea
\det(D^\RR) \det(Z^\RR)^s  = \frac{1}{2^{4n}}(2s-1) (2s)^2(2s+1)^2...(2s+2n-2)^2(2s+2n-1) \det(Z^\RR)^{s-1} 
\eea
}
Just use $\det(M^\RR)=\det(M^\CC)^2$ and previous result.

{\bf Question}
Can similar results be extended to (anti)-symmetric quaternionic matrices ?
The problem that $Z^\CC$, $D^\CC$ will not be (anti)-symmetric
for which Cayley identity is known,
due to very small problem - near diagonal.
E.g. in antisymmetric case we get $
\left( \begin{array}{ccccccccc}
0 & 0\\
0 & 0\end{array} \right)
$ term on the diagonal.
And in symmetric case we get 
$
\left( \begin{array}{ccccccccc}
z_1 & z_2\\
-\bar z_2 &\bar z_1\end{array} \right)
$
which is NOT symmetric, due to $z_2$ and $\bar z_2$.
And the Cayley identities does not seem to be available for such
matrices.

\section{Concluding remarks and further problems \label{SectLast} } 

\subsection{Difficulties with Capelli's "dequanternification" \label{SSectDeq} } 

In our basic motivating paper \cite{AH11} 
"dequaternification" of 1x1 Capelli identity  was obtained.
So 
it is natural
to ask about "dequaternification" for  $n\times n$
matrices.
Of course, we believe that it exists. But it seems our approach
should be somehow modified. Because of the two reasons below.

1) Key thing in our approach is factorization of the
determinants $det(M^\RR)=det(M)det(\bar M)$.
In quaternionic commutative setup 
$\det(M^\RR)= \det(M^\CC)^2$, as it was discussed above.
However as result from \cite{AH11} 
shows in non-commutative setup it is not true that
 $\det(M^\RR)= (something)^2$,
 because his result states  
$\det(M^\RR)=  
\det(X^\RR)^2 \det(D^\RR)^2$ and hence $ \ne
(\det(X^\RR) \det(D^\RR))^2$ due to non-commutativity.

So factorization should be somehow modified.

2) It was important for us that $Z$ and $\bar D$ commute.
In quaternionic case we do not see analogues of this happy circumstance.

So we need to modify somehow the approach to work with quaternionic matrices.

\subsection{Further developments } 

To conclude we have extended the Capelli identities
to "decomplexified" matrices.

The further problems might be the following.

1) Consider "dequaternification".

2) Try to extend linear algebra theorems (Cramer rule, Cayley-Hamilton theorem, ...) for these "decomplexified" and "dequanternified" matrices.
Similar to what have been done in \cite{CFR09} for Manin matrices.

3) Try to extend results about permanents and imananats 
to these cases.

4) Consider super-q analogs of the present results.
(Capelli identities themselves should be 
extendable (and partly already extended) to elliptic and other Lax matrices (see discussion and references in \cite{CF07})).

5) Further examples of  Cayley identities from \cite{CSS11}, \cite{CW81}
might be explored for  "decomplexification" and "dequaternionification".

6) Consider matrices $E^{\RR}$ and its submatrix $S\in E$ which is NOT
a decomplexification of the some submatrix in $E$ - can one extend any
results to such matrix ? 
This is important from general perspective - if matrix $M$ is "good" matrix
(meaning that it has determinant and many linear algebra theorems can be extended
to it) does it mean that any submatrix of $M$ is "good" ? 

7) Our big desire is  to find a Grassman algebra proof for the 
Capelli identity in anti-symmetric case (Howe-Umeda Kostant-Sahi) (see remark \pref{HUKSisHardRemark}),
because we think that such a proof will immediately lead to some generalizations
and hopefully  will extend our understanding of non-commutative linear algebra in general.
It might be that the idea of the present paper "local calculation,
global cancellation" would be useful to attack this problem.

We also mentioned some problems in section \pref{WickIntrSect}.

However ideal solution would be to develop the general formalism
of non-commutative matrices which will include all these possible directions
as a particular case. We hope to make some small step in this direction
in our forthcoming work.


\section{Appendix. Cayley identity via Harish-Chandra's radial part calculation.\label{AppendSect} }
Let us mention the very short proof of the Cayley identity by means of the powerful result - Harish-Chandra's
radial part calculation for certain invariant differential operators. This appendix is independent of the previous text.
We put it here mainly for pedagogical reasons - it seems to us a nice exercise:
short proof of a concrete result by means of general technique which recent days become
quite famous due to P. Etingof and V. Ginzburg \cite{EG00} who
put it into more general context and related to various topics
such as Calogero-Moser integrable hamiltonian system, Cherednik algebras, hamiltonian reduction, n!-conjecture,  etc. 

We were not able to locate such a proof in a literature, so it might be new. 
Similar idea: the proof of another Cayley type identity via the radial part calculation
has been used in a paper by B. Rubin \cite{Ru06}. 
In that paper more recent and complicated identity discussed in \cite{Kh05} 
 has been proved. It would be nice to have proofs of symmetric and antisymmetric analogues
 of the Cayley's identity by radial part calculation and understand what are 
 proper analogues of the Harish-Chandra's isomorphism in these cases.

Let us again refer to a remarkable paper \cite{CSS11} for everything concerning Cayley identity 
(various proofs of old and new results, 
extensive reference list, history account, and relations to the Capelli identity).

\subsection{Reminder on the  Cayley identity and Harish-Chandra's result.}
The Cayley identity has been already discussed, but let us remind it again to make this appendix
independent from the other parts.




Consider the polynomial algebra $\CC[x_{ij}]$, and the matrices:
\bea
X= \left( \begin{array}{ccc}
x_{11} & ... & x_{1n} \\
... & ... & ... \\
x_{n1} & ... & x_{nn} \\
\end{array} \right), ~~~
\partial = \left( \begin{array}{ccc}
\frac{\partial}{\partial x_{11}}  & ... & \frac{\partial}{\partial x_{n1}} \\
... & ... & ... \\
\frac{\partial}{\partial x_{1n}}  & ... & \frac{\partial}{\partial x_{nn}} \\
\end{array} \right). ~~~
\eea
The following beautiful identity is usually attributed to A. Cayley
(we refer to 
\href{http://arxiv.org/abs/0809.3516}{(S. Caracciolo, A. Sportiello, A. Sokal)}
 \cite{CSS11} for the comprehensive treatment of the subject
 and the historical notes):

\bea
\det( \partial ) \det( X )^s =  
s(s+1)...(s+n-1) \det( X )^{s-1}.
\eea

We aim to give a conceptually and techinically few lines short proof of this identity, the price we pay is a refernce to non-trivial and somewhat mystrerious result (attributed to?) by Harish-Chandra on the radial part calculation of the invariant differential operators.
(Main source for it is the famous paper by  
\href{http://arxiv.org/abs/math/0011114}{(P. Etingof, V. Ginzburg)}
\cite{EG00}).




{\bf Fact (Toy model for the Harish-Chandra's result
 ).} Restriction of the 3-dimensional Laplace operator
$\Delta = \frac{\partial^2}{\partial x^2} + \frac{\partial^2}{\partial y^2} +\frac{\partial^2}{\partial z^2}$
 on the rotation invariant functions $f(r)$ equals to the radial Laplacian (i.e. $\frac{\partial^2}{\partial r^2}$) conjugated by "$r$" (i.e. square root of  the volume of the $SO(3)$ orbit):
\bea
\Delta f(r) =  (\frac{\partial^2}{\partial r^2} + \frac{2}{r} \frac{\partial}{\partial r} ) f(r) = (\frac{1}{r} \frac{\partial^2}{\partial r^2}  r) f(r).
\eea

{\bf Mystery.} The fact holds only in  dimension 3. Even for $\RR^2$ it is not true.

{\Rem ~} $\RR^3$ is isomorphic to Lie algebra $so(3)$ with respect to metrics and  $SO(3)$-actions. Harish-Chandra generalized the fact above to all semisimple Lie algebras.

{\vskip .5 cm}

{\bf {\Large Theorem.} Harish-Chandra's radial part calculation for invariant differential operators on Lie algebra $gl_n$}. (See \cite{Et06}  Proposition 4.5 page 27).

Consider Lie algebra $gl_n$, Lie group $GL(n)$ acts on it by conjugation $g\mapsto GgG^{-1}$. 
For an invariant polynomial $f(x_{ij})$ on $gl_n$ denote by $f_{res}(\lambda_1,...,\lambda_n)$ its restriction on the space (slice)  of diagonal matrices.  Respectively for 
an invariant differential operator   {\bf with constant } coefficients  
$D(\partial_{x_{ij}})$ denote by $D_{res}(\partial_{\lambda_1},...,\partial_{\lambda_1})$ similar restriction for differential operators.
Denote by $V(\lambda_1,...,\lambda_n) = \prod_{i<j} (\lambda_i-\lambda_j)$ the
Vandermonde determinant.
Then the following is true:

\bea
(D f)_{res} = \frac{1}{ V(\lambda_1,...,\lambda_n) } D_{res} V(\lambda_1,...,\lambda_n)   f_{res}.
\eea

{\Rem ~} Actually Harish-Chandra's result is true for any Lie algebra, respectively Vandermonde will be
substituted by the Weyl denominanotor (\cite{EG00}).

{\Rem ~} Differential operators with {\bf constant} coefficients means  operators linearly spanned by  $C 
\prod_{ij} (\frac{\partial} { \partial x_{ij} })^{k_{ij}} $,
where $C\in \CC$, i.e. $C$ does not depend on $x_{ij}$.
Clearly this subspace of operators is invariant with respect to $GL(n)$ and
 is isomorphic as representation of $GL(n)$ to the space of polynomials in $x_{ij}$,
polynomials can be restricted to diagonal matrices, hence it makes sense to speak about  similar restiction for these operators.

{\Rem ~} Clearly all invariant functions can be obtained as products of traces  $Tr (X^k)$. The restriction $Tr(X^k)$ on the diagonal matrices is clearly $\sum_{i=1...n} \lambda_i^k$.
Similar all invariant differential  operators with {\bf constant} coefficients can be obtained as products of  $Tr (\partial^k)$. 
The restriction $D_{res}$ of  $Tr(\partial^k)$ is clearly $\sum_{i=1...n} (\frac{\partial}{\partial   \lambda_{i} }) ^k$, the point of Harish-Chandra's result that this naive restriction does not give correct action on invariant function, but it differs by quite a manageable expression.

{\Ex ~} 
Consider $gl_2$, take $D= Tr(\partial) = \frac{\partial}{\partial x_{11} } +
\frac{\partial}{\partial x_{22} } $; take $f=Tr(X) =  x_{11} +x_{22} $.

\bea
(D f )_{res} = ( (\frac{\partial}{\partial x_{11} } +
\frac{\partial}{\partial x_{22} } ) ( x_{11} +x_{22} ) )_{res} = 
(1 + 1)_{res} = 2.
\eea

\bea
\frac{1} { V(\lambda_1,  \lambda_2) } ~ D_{res} ~ V(\lambda_1 , \lambda_2) ~  (f )_{res} 
=
\frac{1} { (\lambda_1 - \lambda_2) } ~( \frac{\partial }{\partial \lambda_{1} } 
+ \frac{\partial }{\partial \lambda_{2} }  ) ~ (\lambda_{1}  - \lambda_{2}  ) ~ (\lambda_{1}  + \lambda_{2}  ) 
= \\
=
\frac{1} { (\lambda_1 - \lambda_2) } ~ ( \frac{\partial }{\partial \lambda_{1} } 
+ \frac{\partial }{\partial \lambda_{2} }  )  ~ (\lambda_{1}^2  - \lambda_{2}^2  ) 
=
\frac{1} { (\lambda_1 - \lambda_2) }  (2 \lambda_{1}  - 2 \lambda_{2}  ) 
= 2.
\eea

{\Rem ~}
The appearance of the square root of the volume of the orbit is natural by several reasons, for example because the scalar products on the total space and the  slice
differs by the volume of the orbit, so to get the isomorphism of the Hilbert  space of invariant functions on the total space and the Hilbert space of functions on the slice one should multiply by the square root of the orbit volume. 
Somewhat mysterious is that the result holds exactly in Harish-Chandra's setup while for other cases holds only approximately.

\subsection{Proof of the Cayley identity.}
{\PRF}
Due to Harish-Chandra's result we need to prove the following:
\bea
\label{CayleyReduced}\frac{1} { V(\lambda_1, ... , \lambda_n) } 
(\frac{\partial} {\partial \lambda_1} \frac{\partial} {\partial \lambda_2}...
 \frac{\partial} {\partial \lambda_n} ) 
 { V(\lambda_1, ... , \lambda_n) } 
 (\lambda_1  \lambda_2... \lambda_n )^s = 
s(s+1)(s+2)...(s+n-1)   (\lambda_1  \lambda_2... \lambda_n )^{s-1}.
\eea

Recall that: 
\bea  
 V(\lambda_1, ... , \lambda_n) = \det \left( \begin{array}{cccc}
(\lambda_1)^0 & (\lambda_2)^0 & ... & (\lambda_n)^0  \\
(\lambda_1)^1 & (\lambda_2)^1 & ... & (\lambda_n)^1  \\
... & ... & ... & ... \\
(\lambda_1)^{n-1} & (\lambda_2)^{n-1} & ... & (\lambda_n)^{n-1}  \\
\end{array} \right) 
= 
\sum_{\sigma \in S_n} (-1)^{sgn(\sigma)}  
\prod_{i=1...n} \lambda_{i}^{\sigma(i) -1 } .
\eea

Observe that operator $    
(\frac{\partial} {\partial \lambda_1} \frac{\partial} {\partial \lambda_2}...
 \frac{\partial} {\partial \lambda_n} ) 
$
acting on any monomial $\prod_{i} (\lambda_i)^{k_i} $ produces 
factor $ \prod_i k_i$ and lowers all degrees by one.

Now observe that for  all monomials which enter the expansion of the Vandermonde 
multiplied by $(\lambda_1  \lambda_2... \lambda_n )^{s-1}$
give the same factor 
$ \prod_i k_i=s(s+1)...(s+n-1)$ because  degrees $k_i$ of $\lambda_i$ are permutation of the same
set of numbers $s, (s+1), ..., (s+n-1)$ so the product is always the same.

Taking into account that degrees are lowered by one we conclude that:
\bea
(\frac{\partial} {\partial \lambda_1} \frac{\partial} {\partial \lambda_2}...
 \frac{\partial} {\partial \lambda_n} ) 
 { V(\lambda_1, ... , \lambda_n) } 
 (\lambda_1  \lambda_2... \lambda_n )^s = 
s(s+1)(s+2)...(s+n-1) { V(\lambda_1, ... , \lambda_n) }   (\lambda_1  \lambda_2... \lambda_n )^{s-1}.
\eea
This gives the desired formula \ref{CayleyReduced} just diving by the Vandermonde.

Proof is finished.
\BX

{\bf Acknowledgements.}
The author is grateful to F. Khafizov, D. Shmelkin, S. Slavnov,
R. Bebenin, B. Shulgin without whom this paper would never exist,
for discussions, support and inspiration; to P. Etingof for clarifications, to A. Molev for useful comments. 
I am thankful to  V. Dotsenko,  M. Sapir, D. Speyer, M. Tucker-Simmons for discussions on Mathoverflow
(see footnote on  page \pageref{thankMO}).
Special thanks to my family: wife Almira, daughter Kamila, mother Lily
for their support and understanding.
The author thanks grant RFBR-090100239a.

\end{document}